\documentclass[11pt, reqno]{amsart}

\usepackage{amssymb,amsmath,eufrak}
\author[Florent Benaych-Georges]{Florent Benaych-Georges}
\title[Infinitely divisible distributions for rectangular free convolution]{Infinitely divisible distributions for rectangular free convolution: classification and matricial interpretation}
\date{\today}
\newcommand{\bxpp}{{\scriptscriptstyle\boxplus}}

\newcommand{\Part}{\operatorname{Part}}
\newcommand{\NC}{\operatorname{NC}}
\newcommand{\NCp}{\operatorname{NC'}}

\newcommand{\ds}{\displaystyle}

\newcommand{\tr}{\operatorname{tr}}
\newcommand{\rg}{\operatorname{rg}}

\newcommand{\mf}{\mathfrak}

\newcommand{\Tr}{\operatorname{Tr}}
\newcommand{\ninf}{\underset{n\to\infty}{\longrightarrow}}
\newcommand{\ssi}{if and only if }

\newcommand{\teo}{theorem }
\newcommand{\teov}{theorem, }

\newcommand{\csts}{$\ast$-stable }

\newcommand{\cid}{$\ast$-infinitely divisible distribution }
\newcommand{\cids}{$\ast$-infinitely divisible distributions }

\newcommand{\fids}{$\sst\boxplus$-infinitely divisible distributions }

\newcommand{\cidp}{$\ast$-infinitely divisible distribution. }

\newcommand{\rid}{$\scriptstyle\boxplus_\la$-infinitely divisible distribution }

\newcommand{\E}{\mathbb{E}}

\newcommand{\R}{\mathbb{R}}
\newcommand{\C}{\mathbb{C}}
\newcommand{\n}{\mathbb{N}}

\newcommand{\ud}{\mathrm{d}}

\newcommand{\trv}{Voiculescu transform }

\newcommand{\pro}{probability }

\newcommand{\pmu}{ $\mathbb{P}_{d,d'}^\mu$ }
\newcommand{\pmum}{ \mathbb{P}_{d,d'}^\mu }

\newcommand{\rdmunp}{ $\mathbb{Q}_{d,d'}^{\mu_n}$. }

\newcommand{\dtpd}{$d\scriptstyle\times\ds\!  d'$ }

\newcommand{\f}{\frac}
\newcommand{\ff}{\frac{1}}
\newcommand{\lf}{\left}
\newcommand{\ri}{\right}

\newcommand{\st}{such that }
\newcommand{\la}{\lambda}
\newcommand{\La}{\Lambda}
\newcommand{\tii}{\scriptstyle\times\ds\!}

\newcommand{\vfi}{\varphi}
\newcommand{\ste}{\, ;\, }
\newcommand{\mc}{\mathcal }

\newcommand{\eps}{\varepsilon}
\newcommand{\arc}{{\scriptscriptstyle\boxplus_{\la}}}
\newcommand{\arco}{{\scriptscriptstyle\boxplus_{0}}}
\newcommand{\gab}{\Delta_{\alpha,\beta}}
\newcommand{\sst}{\scriptstyle}

\newcommand{\bxp}{{\scriptscriptstyle\boxplus}}

\newcommand{\D}{\mc{D}}

\numberwithin{equation}{section}
\newtheorem{Th}{Theorem}[section]
\newtheorem{propo}[Th]{Proposition} 
 
\newtheorem{lem}[Th]{Lemma}

\newtheorem{rmq}[Th]{Remark}
\newtheorem{cor}[Th]{Corollary}
\newtheorem{Def}[Th]{Definition}

\newenvironment{pr}{\noindent {\bf Proof. }}{\ \ \ $\square$}
\newenvironment{prth}{\noindent {\bf Proof of the theorem. }}{\ \ \ $\square$}

\long\def\symbolfootnote[#1]#2{\begingroup
\def\thefootnote{\fnsymbol{footnote}}\footnote[#1]{#2}\endgroup}

\begin{document}
\maketitle
\symbolfootnote[0]{{\it MSC 2000 subject classifications.}  15A52, 
46L54, 
60E07, 
60F05} 
 
\symbolfootnote[0]{{\it Key words.} random matrices, free probability, free convolution, Marchenko-Pastur distribution, infinitely divisible distributions} 

\begin{abstract}In a previous paper (\cite{fbg.AOP.rect}), we defined the rectangular free convolution $\arc$. Here, we investigate  the related notion of infinite divisibility, which happens to be closely related the classical infinite divisibility: there exists a bijection between the set of classical symmetric infinitely divisible distributions and the set of  $\arc$-infinitely divisible distributions, which preserves limit theorems. We give an interpretation of this correspondence in terms of random matrices: we construct distributions on sets of complex rectangular matrices which give rise to random matrices with singular laws going from the symmetric classical infinitely divisible distributions to their  $\arc$-infinitely divisible correspondents when the dimensions go from one to infinity in a ratio $\la$.\end{abstract}\tableofcontents

\section*{Introduction} In   a previous paper (\cite{fbg.AOP.rect}), we modeled the asymptotic behavior of rectangular random matrices with freeness with amalgamation.
Therefore we defined, for each $\la\in [0,1]$, the rectangular free convolution with ratio $\la$, denoted by $\arc$. It is a binary operation on the set of symmetric \pro measures on the real line defined in the following way. Let us call the singular law of a matrix $M$ the uniform law on its singular values, i.e. on the spectrum of its absolute value $|M|=(MM^*)^{1/\!2}$. 
Consider $\mu,\nu$  symmetric \pro measures on the real line, consider two sequences $q_1(n),q_2(n)$ of  integers tending to $+\infty$ \st $$\f{q_1(n)}{q_2(n)} \ninf \la,$$ and consider, for each $n$, $M(n), N(n)$ independent $q_1(n)\tii q_2(n)$ random matrices, one of them being biunitarily invariant (i.e. having a distribution invariant under the left and right actions of the unitary groups) \st the symmetrization of the singular law of $M(n)$ (resp. of $N(n)$) converges weakly in \pro to $\mu$ (resp. $\nu$).  Then the symmetrization of the singular law of $M(n)+N(n)$ converges weakly in \pro to a \pro measure which depends only on $\mu$, $\nu$, and $\la$, denoted by $\mu\arc\nu$, and called the {\it rectangular free convolution with ratio $\la$} of $\mu$ and $\nu$.

In the present paper, we study the notion of infinite divisibility for $\arc$, which leads to a L\'evy-Kinchine formula for the rectangular $R$-transform (whose definition we shall recall in section \ref{archi.nancy+antropo.funeraire}): a symmetric \pro measure $\mu$ is  $\arc$-infinitely divisible \ssi there exists a positive finite symmetric measure $G$ (called its {\it L\'evy measure}) \st the rectangular $R$-transform with ratio $\la$ of $\mu$ is given by the formula: $$C_\mu(z)=z\int_\R\f{1+t^2}{1-zt^2}\ud G(t).$$  

Therefore we can define a bijection $\La_\la$ between the set of classical symmetric infinitely divisible distributions and the set of $\arc$-infinitely divisible   distributions:  $\La_\la$ maps a symmetric $*$-infinitely divisible distribution to the  $\arc$-infinitely divisible   distribution with the same L\'evy measure. This bijection happens, like the one of Bercovici and Pata (\cite{appenice}) between $*$- and $\bxp$-infinitely divisible distributions, to have deep properties. It is a semi-group morphism: $$\La_\la(\mu*\nu)=\La_\la(\mu)\arc\La_\la(\nu),$$ and it preserves limit theorems: for all sequences $(\mu_n)$ of symmetric distributions and $(k_n)$ of positive integers tending to infinity, we have, for all \pro measures $\mu$, $$\mu_n^{*k_n}\ninf \mu\iff \mu_n^{\arc k_n}\ninf\La_\la(\mu).$$ $\La_\la$ will be called the rectangular Bercovici-Pata bijection with ratio $\la$.

In section \ref{stzn.21.09.06}, we give examples of $\arc$-infinitely divisible distributions. 
First, in section \ref{6.2.05.1} we give the density of the image, by the bijection $\La_\la$, of the standard Gaussian distribution. 
An interesting interpretation of this result is made in  a forthcoming paper (\cite{fbg.free.amalg}) where we construct analogues of Voiculescu's free entropy and free Fisher information for operators between different Hilbert spaces, and where the maximum of entropy and the minimum of Fisher information are realized for operators the absolute value of which has this (symmetrized) distribution. Another consequence of this result is a new formula for the moments of the Marchenko-Pastur distribution (which is closely related to this distribution): for all $a>0$, for all $n\geq 1$, the $n$-th moment of the Marchenko-Pastur distribution with parameter $a$ (see \cite{hiai} p. 65) is equal to $\sum_{\pi}a^{o(\pi)}$, where the sum is taken over all noncrossing pairings of $[2n]$, and where $o(\pi)$ is the number of blocks of a partition $\pi$ the first element of which is odd. 
Then, in section \ref{cauchysection.22.09.06}, we give the densities of the  images, by the above mentioned bijection, of the symmetric Cauchy laws. Cauchy laws are well known to be invariant under  many transformations, but we are going to see that unless $\la=1$, they are not invariant under this bijection. At last, in section \ref{bon.vieux.dede.20.09.06}, we characterize the images, by the same bijection, of the symmetric Poisson distributions. When $\la=0$, we have a formula for the density. 

In section \ref{broderie+cafe=champion}, we   shall construct a matricial model  for the $\arc$-infinitely divisible laws and present in a maybe more palpable way the Bercovici-Pata bijection with ratio $\la$ (whereas the proofs of the other sections rely on integral transforms and complex analysis):  we are going to construct, in the same way as in \cite{fbg.ID} and in \cite{cabduv.ID}, for each $d,d'\geq 1$, for each symmetric \cid $\mu$, an infinitely divisible  distribution \pmu on the set of \dtpd complex matrices \st for all $\mu,\nu$, $\mathbb{P}_{d,d'}^\mu*\mathbb{P}_{d,d'}^\nu=\mathbb{P}_{d,d'}^{\mu *\nu}$ and \st the symmetrization of the singular law  of  $M$ (with $M$ random matrix distributed according to $\mathbb{P}_{d,d'}^\mu$)  goes from $\mu$ to its image by the rectangular Bercovici-Pata bijection with ratio $\la$ when $d,d'\to\infty , \f{d}{d'}\to\la$.

In the last section, we shall give a representation of the image of the symmetric Poisson distribution by  the rectangular Bercovici-Pata bijection with ratio $\la$ as the distribution of the absolute value of sums of rank-one matrices.

{\bf Acknowledgments.} We would like to thank Philippe Biane, our advisor, for useful discussions. Also, we would like to thank C\'ecile Martineau for her contribution to the english version of this paper.

\section{Preliminaries}\label{archi.nancy+antropo.funeraire}
Until the end of section \ref{stzn.21.09.06},  $\la$ is a fixed number of $[0,1]$. 
\subsection{General introduction to the rectangular $R$-transform with ratio $\la\in [0,1]$}
In this section, we shall recall definitions and basic results from  \cite{fbg.AOP.rect} about the rectangular $R$-transform $C_\mu$ of a symmetric \pro measure $\mu$.

Let us denote by $z\mapsto z^{1/\!2}$ (resp. $z\mapsto \sqrt{z}$) the analytic version of the square root on the complement of the real non positive (resp. non negative) half line \st $1^{1/\!2}=1$ (resp. $\sqrt{-1}=i$). On the set of non null complex numbers, we will use   the argument function which takes values in $[0,2\pi)$.  Let us define the analytic function  on a neighborhood of zero $U(z)=  \f{-\la-1+\lf[(\la+1)^2+4\la z\ri]^{1/2}}{2\la} $ (when $\la=0$, $U(z)=z$).  Then one can summarize the different steps of the construction of the rectangular $R$-transform with ratio $\la$ in the following chain$$
\begin{array}{l}\ds\underset{\substack{{\textrm{sym. prob.}}\\ {\textrm{measure}}}}{\mu}\,\,\longrightarrow \,\,\underset{{\textrm{Cauchy transform}}}{G_\mu(z)=\int\f{\ud \mu(t)}{z-t}}\,\,\longrightarrow\,\, H_\mu(z)=\la G_\mu\lf(\ff{\sqrt{z}}\ri)^2+(1-\la)\sqrt{z}G_\mu\lf(\ff{\sqrt{z}}\ri)\,\,\longrightarrow\,\,\\ \ds  \underset{\textrm{rect. $R$-transf. with ratio $\la$}}{C_\mu(z)=U\lf( \f{z}{H_\mu^{-1}(z)}-1\ri),}\end{array}$$where $H_\mu^{-1}$ is the inverse (for composition) of $H_\mu$. Proposition \ref{Delphine=princessecornichon} and theorem \ref{passiontroisiemeage} bellow, which have first been established in \cite{fbg.AOP.rect}, prove that such an inverse exists, give its domain, and prove that for any tight set $A$ of symmetric \pro measures,  the properties and the domains of the functions $H_\mu^{-1}$ ($\mu\in A$) are ``uniform''. 
\begin{propo}\label{Delphine=princessecornichon} Let $A$ be a set of symmetric \pro measures on the real line. Then the following assertions are equivalent

\begin{itemize}
\item[(i)] A is tight,
\item[(ii)] for every $0<\theta < \pi$,  $\underset{\substack{z\to 0\\ \lf|\arg z -\pi\ri|<\theta}}{\lim}\ff{z}H_\mu(z)=1$ uniformly in $\mu\in A$,
\item[(iii)] $\underset{\substack{x\to 0\\ x\in (-\infty,0)}}{\lim}\ff{x} H_\mu(x) =1$  uniformly in $\mu\in A$.
\end{itemize}
\end{propo}

Define, for $\alpha \in (0, \pi)$, $\beta > 0$, $\Delta_{\alpha,\beta}$ to be the set of complex numbers $z$ \st $|\arg z-\pi|< \alpha $ and $|z|<\beta$.

Let $\mc{H}$ be the set of functions $f$ which are analytic in a domain $\D_f$ \st for all  $\alpha \in (0, \pi)$, there exists $\beta$ positive \st $$\Delta_{\alpha,\beta}\subset\D_f.$$ A family $(f_a)_{a\in A}$ of functions of $\mc{H}$ is said to be {\it uniform} if for all  $\alpha \in (0, \pi)$,  there exists $\beta$ positive \st $$\forall a\in A,\quad \Delta_{\alpha,\beta}\subset\D_{f_a}.$$

\begin{Th}\label{passiontroisiemeage}
Let $(H_a)_{a\in A}$ be a uniform family of functions of $\mc{H}$ \st for every $\alpha\in (0, \pi)$,  $$\underset{\substack{z\to 0\\ |\arg z-\pi|< \alpha }}{\lim}\f{H_a(z)}{z}=1 \textrm{ uniformly in $a\in A$.}$$

Then there exists a uniform family $(F_a)_{a\in A}$ of functions  of $\mc{H}$ \st for every $\alpha\in (0, \pi)$,  $$\underset{\substack{z\to0\\ |\arg z-\pi|< \alpha }}{\lim}\f{F_a(z)}{z}=1\textrm{ uniformly in $a\in A$,}$$ and there exists $\beta$ positive \st $$\forall a\in A,\quad H_a\circ F_a=F_a\circ H_a=I_d\textrm{ on $\Delta_{\alpha,\beta}$.}$$ 

Moreover, the family $(F_a)_{a\in A}$ is unique in the following sense: if a family $(\tilde{F}_a)_{a\in A}$  of functions  of $\mc{H}$ satisfies the same conditions, then for all  $\alpha \in (0, \pi)$, there exists $\beta$ positive \st $$\forall a\in A,\quad F_a=\tilde{F}_a \textrm{ on }\Delta_{\alpha,\beta}.$$
\end{Th}

Using the theory of cumulants in operator-valued free \pro theory, we prove  (\cite{fbg.AOP.rect}) the {\it additivity of rectangular $R$-transform}: \begin{Th}\label{additivity.21.09.06}For all $\mu,\nu$, symmetric \pro measures, we have$$C_{\mu\arc\nu}=C_\mu+C_\nu.$$\end{Th}

Note that it is also proved in \cite{fbg.AOP.rect} that for all $\la\in [0,1]$, the rectangular $R$-transform with ratio $\la$ is injective. The following remark gives a practical way to derive  any symmetric \pro measure $\mu$ on the real line  from $C_\mu$.

\begin{rmq}[How to compute $\mu$ when we know $C_\mu$ ?]\label{21.09.06.1} Let us define the function $T(X)=(\la X+1)(X+1)$, \st $T(U(z))=z+1$. We have $z/H_\mu^{-1}(z)=T(C_\mu(z))$, for $z\in \C\backslash \R^+$ small enough. From this, we can compute $H_\mu(z)$ for $z\in \C\backslash \R^+$ small enough. Then we can use the equation, for  $z\in \C\backslash \R^+$, $$\ff{z}H_\mu(z)=\la\lf(\ff{\sqrt{z}}G_\mu\lf(\ff{\sqrt{z}}\ri)\ri)^2+(1-\la)\ff{\sqrt{z}}G_\mu\lf(\ff{\sqrt{z}}\ri).$$Moreover, when $z\in \C\backslash \R^+$ is small enough, $1/\sqrt{z}$ is large and in $\C^-$ (the set of complex numbers with negative imaginary part), so $\ff{\sqrt{z}}G_\mu\lf(\ff{\sqrt{z}}\ri)$ is closed to $1$. $\ff{z}H_\mu(z)$ is also closed to $1$, and for $h,g$ complex numbers closed to $1$, $$h=\la g^2+(1-\la)g\Leftrightarrow g=V(h),\textrm{ with }V(z)=\f{\la-1+((\la-1)^2+4\la z)^\ff{2}}{2\la}=U(z-1)+1.$$ 
So one has, for $z\in  \C\backslash \R^+$ small enough,  $$\ff{\sqrt{z}}G_\mu\lf(\ff{\sqrt{z}}\ri)=V\lf(\f{H_\mu(z)}{z}\ri).$$ \end{rmq}

We shall mention here two other  results, proved in  \cite{fbg.AOP.rect}. The second of them allows us to claim that $\arc $ is continuous with respect to weak convergence.
 \begin{lem}[Tightness and rectangular $R$-transform]\label{18.11.03.2} Let $A$ be a set of symmetric \pro measures. Then we have equivalence between :

\begin{itemize}
\item[(i)] $A$ is tight,
\item[(ii)] for any $0<\alpha <\pi$,  $\underset{\substack{z\to 0\\ \lf|\arg z -\pi\ri|<\alpha}}{\lim}C_\mu(z)=0$ uniformly in $\mu\in A$,
\item[(iii)] $\underset{\substack{x\to 0\\ x\in (-\infty,0)}}{\lim}C_\mu(x)=0$ uniformly in $\mu\in A$.
\end{itemize} 
\end{lem}

\begin{Th}[Paul L\'evy's \teo for rectangular $R$-transform]\label{tesprof?}  
Let $(\mu_n)$ be a sequence of symmetric \pro measures. Then we have equivalence between:
\begin{itemize}
\item[(i)]  $(\mu_n)$ converges weakly to a symmetric \pro measure;
\item[(ii)] there exists $\alpha,\beta$ \st \begin{itemize}\item[(a)]   $\underset{\substack{z\to 0\\ |\arg z-\pi|<\alpha}}{\lim}C_{\mu_n}(z)=0$ uniformly in $n$, \item[(b)] the sequence $(C_{\mu_n})$ converges uniformly on every compact set of $\Delta_{\alpha,\beta}$ when $n\to\infty$;\end{itemize} 
\item[(iii)] \begin{itemize}\item[(a)] $\underset{\substack{x\to 0\\ x\in (-\infty,0)}}{\lim}C_{\mu_n}(x)=0$ uniformly in $n$,  \item[(b)] there exists $\beta> 0$ \st the sequence $(C_{\mu_n})$ converges pointwise  on $[-\beta, 0)$ when $n\to\infty$. \end{itemize} 
\end{itemize} Moreover, in this case, denoting by $\mu$ the weak limit of $(\mu_n)$, for every $\alpha$, there exists $\beta$ \st the sequence $(C_{\mu_n})$ converges uniformly to $C_\mu$ on every compact set of $\Delta_{\alpha,\beta}$ when $n\to\infty$. 
\end{Th}

\subsection{The particular cases $\la=0$ and $\la=1$}\label{20.09.06.1}The results of this section are proved in \cite{fbg.AOP.rect}. 

\subsubsection{Rectangular free convolution}\label{foule.romaine.19.09.06} For $\mu,\nu$ symmetric \pro measures on the real line, the rectangular free convolution with ratio $1$ of $\mu$ and $\nu$ is their free convolution (as defined in \cite{defconv}), and their rectangular free convolution with ratio $0$ is the unique symmetric \pro measure on the real line whose push-forward by the function $t\to t^2$ is the free convolution of their push-forwards by the same function.

\subsubsection{Rectangular $R$-transform} The rectangular $R$-transform with ratio $1$ (resp. $0$), for a symmetric distribution $\mu$,  is linked to  the \trv $\vfi_\mu$ of $\mu$  by the relation $C_\mu(z)=\sqrt{z}\vfi_\mu(1/\sqrt{z})$ (resp. $C_\mu(z)=z\vfi_\rho(1/z)$, where $\rho$ is the push-forward of $\mu$ by the function $t\to t^2$) (see paragraph 5 of \cite{defconv} for the construction of the Voiculescu transform).

\section{L\'evy-Kinchine Formula for $\arc$-infinitely divisible distributions}
$\arc$-infinitely divisible distributions are defined in the same way as $*$- and ${\sst \boxplus}$-infinitely divisible distributions: 
\begin{Def} A symmetric \pro measure $\nu$ is said to be  $\arc$-infinitely divisible if for each $n\in \n^*$, there exists a symmetric distribution $\nu_n$ \st $\nu_n^{\arc n}=\nu$. 
\end{Def}
 As for $*$- and ${\scriptstyle\boxplus}$-, we have the following characterization of $\arc$-infinite divisibility. 

\begin{Th} Let $\nu$ be a symmetric distribution. Then $\nu$ is $\arc$-infinitely divisible \ssi there exists a sequence $(\nu_n)$ of symmetric \pro measures  \st $\nu_n^{\arc n}$ converges weakly to $\nu$. 
\end{Th}

\begin{pr}If $\nu$ is  $\arc$-infinitely divisible, it is clear. Assume the existence of a sequence $(\nu_n)$ \st $\nu_n^{\arc n}$ converges weakly to $\nu$. Consider $k\geq 1$. Let us show that there exists a symmetric \pro measure $\sigma$ \st $\sigma^{\arc k}=\nu$. We have $\ds\lim_{\substack{y\to 0\\ y<0}}nC_{\nu_n}(y)=0$ uniformly in $n$, so $\ds\lim_{\substack{y\to 0\\ y<0}}nC_{\nu_{k n}}(y)=0$ uniformly in $n$. So by lemma \ref{18.11.03.2}, the sequence $\lf(\nu_{k n}^{\arc n}\ri)$ is tight. If the symmetric distribution $\sigma$ is the limit of one of its subsequences, we have $$\ds\sigma^{\arc k}=\lim_{n\to \infty} \lf(\nu_{k n}^{\arc n}\ri)^{\arc k}=\lim_{n\to \infty}  \nu_{k n}^{\arc n k}=\nu . $$
\end{pr}
\begin{cor} The set of  $\arc$-infinitely divisible distributions is closed under weak convergence.\end{cor}
\begin{pr}If a sequence $(\mu_n)$ of $\arc$-infinitely divisible distributions converges weakly to a distribution $\mu$, then if for every $n$, $\nu_n^{\arc n}=\mu_n$, the sequence $(\nu_n^{\arc n})$ converges weakly to $\mu$.\end{pr}
\\
\\

To prove the L\'evy-Kinchine formula for $\arc$-infinitely divisible distributions, we need the following lemma, which is the analogue of propositions 2.6 and 2.7 of \cite{appenice}. Until the end of the paper, for $f,g$ functions defined on a domain whose closure contains an element $x_0$,  "$f(x)\sim g(x)$ in the neighborhood of $x_0$ (or as $x$ tends to $x_0$)" will mean that $f(x)/g(x)$ tends to $1$ as $x$ tends to $x_0$.



\begin{lem}\label{27.11.03.3}
Let $(\nu_n)$ be a sequence of symmetric \pro measures that converges weakly to $\delta_0$. Consider $\alpha \in (0,\pi)$. Then there exists $\beta >0$ \st on 
$\Delta_{\alpha,\beta}$, \[C_{ \nu_n}(z)=\lf(\ff{\sqrt{z}}G_{ \nu_n}\lf(\ff{\sqrt{z}}\ri)-1\ri)(1+v_n(z))\]where the functions $v_n$, defined on $\Delta_{\alpha,\beta}$,  are \st
\begin{itemize}\item[(i)] $\forall n, \forall z, |v_n(z)|\leq 1/2$ and $\underset{z\to 0}{\lim}v_n(z)= 0$ uniformly in $n$,
\item[(ii)] $\forall z, \underset{n\to\infty}{\lim}v_n(z)= 0$.
\end{itemize}
\end{lem}

\begin{pr} First, note that unless $\nu_n=\delta_0$ (in which case $v_n=0$ is suitable), for all $z$, $\ff{\sqrt{z}}G_{ \nu_n}\lf(\ff{\sqrt{z}}\ri)-1\neq 0$, so there is a function $v_n$ on the domain of $C_{\nu_n}$ \st \begin{equation}\label{on.a.beach.ashcroft.09.06}C_{ \nu_n}(z)=\lf(\ff{\sqrt{z}}G_{ \nu_n}\lf(\ff{\sqrt{z}}\ri)-1\ri)(1+v_n(z)).\end{equation} To prove $(i)$, we will only use the tightness of $\{\nu_n\ste n\in\n\}$. It suffices to show that
\[\underset{\substack{z\to 0\\ |\arg z-\pi|<\alpha }}{\lim}\f{C_{ \nu_n}(z)}{\ff{\sqrt{z}}G_{ \nu_n}\lf(\ff{\sqrt{z}}\ri)-1}=1 {\textrm{ uniformly in $n$.}}\]
 We have, by the paragraph following Proposition 5.1 in \cite{defconv}, \[\underset{\substack{z\to 0\\ |\arg z-\pi|<\alpha }}{\lim}\ff{\sqrt{z}}G_{ \nu_n}\lf(\ff{\sqrt{z}}\ri)=1 {\textrm{ uniformly in $n$,}}\] and when a complex number $t$ tends to $1$, $t-1\sim(\la t^2+(1-\la )t-1)/(\la +1)$, so it suffices to show that  \[\underset{\substack{z\to 0\\ |\arg z-\pi|<\alpha }}{\lim}\f{(\la +1)C_{ \nu_n}(z)}{\f{H_{ \nu_n}(z)}{z}-1}=1 {\textrm{ uniformly in $n$.}}\] We have $$\f{H_{ \nu_n}(z)}{z}-1=\f{H_{ \nu_n}(z)}{z}\lf(1-\f{z}{H_{ \nu_n}(z)} \ri),$$
and we know, by proposition \ref{Delphine=princessecornichon}, that\begin{equation}\label{Sissi.a.du.boulot1}\underset{\substack{z\to 0\\ |\arg z-\pi|<\alpha }}{\lim}\f{H_{\nu_n}(z)}{z}=1 {\textrm{ uniformly in $n$.}}\end{equation} 
So it suffices to show that  
\[\underset{\substack{z\to 0\\ |\arg z-\pi|<\alpha }}{\lim}\f{(\la +1)C_{ \nu_n}(z)}{1-\f{z}{H_{ \nu_n}(z)}}-1=0 {\textrm{ uniformly in $n$.}}\]  
We know, by proposition  \ref{Delphine=princessecornichon} and by  \teo \ref{passiontroisiemeage},  that  
\begin{equation}\label{Sissi.a.du.boulot2}\underset{\substack{z\to 0\\ |\arg z-\pi|<\alpha }}{\lim}\f{H_{\nu_n}^{-1}(z)}{z}=1 {\textrm{ uniformly in $n$,}}\end{equation} 
and the equivalent of $U(x)$ in a neighborhood of zero is  $\f{x}{\la+1}$. So, since  $\ds C_{\nu_n}(z)=U\lf(\f{z}{H_{ \nu_n}^{-1}(z)}-1\ri)$,  it suffices to show that  
\[
\underset{
\substack{z\to 0\\ |\arg z-\pi|<\alpha }
}{\lim}\f{\f{z}{H_{ \nu_n}^{-1}(z)}-1
}{1-\f{z}{H_{ \nu_n}(z)}}=1 
{\textrm{ uniformly in $n$.}}\] 
Choose $\alpha'\in (\alpha, \pi)$. By theorem \ref{passiontroisiemeage}, there exists  $\beta_1 >0$ \st for all $n$, $H_{\nu_n}^{-1}$ is defined on $\Delta_{\alpha',2\beta_1}$, and $$H_{\nu_n}(\Delta_{\alpha,\beta_1})\cup H_{\nu_n}^{-1}(\Delta_{\alpha,\beta_1})\subset \Delta_{\alpha',2\beta_1}.$$We have, for $z\in\Delta_{\alpha,\beta_1}$,   \begin{eqnarray*}\f{\f{z}{H_{ \nu_n}^{-1}(z)}-1
}{1-\f{z}{H_{ \nu_n}(z)}}-1&=&\f{\f{zH_{ \nu_n}(z)}{H_{ \nu_n}^{-1}(z)}-H_{ \nu_n}(z)
}{H_{ \nu_n}(z)-z}-1\\ 
&=&\ff{H_{ \nu_n}(z)-z}\int_{[z,H_{ \nu_n}(z)]}f_{n,z}'(\xi)\ud \xi, \end{eqnarray*} where $f_{n,z}$ is the function defined by $$f_{n,z}(\xi)=\f{H_{ \nu_n}^{-1}(\xi)H_{ \nu_n}(z)}{H_{ \nu_n}^{-1}(z)}-\xi.$$ 
But 
 the lemma 2.4 of \cite{appenice} states that \[\underset{\substack{z\to 0\\ |\arg z-\pi|<\alpha'}}{\lim}(H_{\nu_n}^{-1})'(z)=1 {\textrm{ uniformly in $n$.}}\] Hence, using also (\ref{Sissi.a.du.boulot1}) and (\ref{Sissi.a.du.boulot2}) (which stay true if $\alpha$ is replaced by $\alpha'$), we  have   \[\ds\underset{
\substack{z\to 0\\ |\arg z-\pi|<\alpha }
}{\lim}\sup\{\lf|f_{n,z}'(\xi )\ri|\ste \xi\in [z, H_{\nu_n}(z)]\}=0 {\textrm{ uniformly in $n$.}}\] so
\[
\underset{
\substack{z\to 0\\ |\arg z-\pi|<\alpha}
}{\lim}\f{\f{z}{H_{ \nu_n}^{-1}(z)}-1
}{1-\f{z}{H_{ \nu_n}(z)}}=1 
{\textrm{ uniformly in $n$.}}\]  So we know that the sequence $(v_n)$ of functions satisfying (\ref{on.a.beach.ashcroft.09.06}) satisfies  $\underset{z\to\infty}{\lim}v_n(z)= 0$ uniformly in $n$. Hence we can choose $\beta_2$ \st $\forall n, \forall z\in\Delta_{\alpha,\beta_2}, |v_n(z)|\leq 1/2$, and (i) is satisfied.

Let us now prove $(ii)$.  First, note that since $\nu_n\to\delta_0$, the sequence $(G_{\nu_n})$ converges uniformly to $G_{\delta_0}: z\mapsto 1/z$ on every compact of the upper half plane (see \cite{akhi} or section 3.1 of \cite{hiai}), so, as in the proof of $(i)$, it suffices to show that\[\ds\lim_{n\to\infty}\f{(\la +1)C_{ \nu_n}(z)}{\f{H_{ \nu_n}(z)}{z}-1}=1.\] The convergence of $\nu_n$ to $\delta_0$ implies too that  $(H_{\nu_n})$ converges to $H_{\delta_0}: z\mapsto z$ . So, since $$\f{H_{ \nu_n}(z)}{z}-1=\f{H_{ \nu_n}(z)}{z}\lf(1-\f{z}{H_{ \nu_n}(z)} \ri),$$ it suffices to prove that \[\lim_{n\to\infty}\f{(\la +1)C_{\nu_n}(z)}{1-\f{z}{H_{ \nu_n}(z)}}=1.\]Furthermore, by \teo \ref{tesprof?}, there exists $\beta_3\leq  \beta_2$ \st the sequence $C_{\nu_n}$ converges uniformly to $C_{\delta_0}=0$ on every compact of $\Delta_{\alpha,\beta_3}$. 
So $z/H_{\nu_n}^{-1}(z)=(\la C_{\nu_n}(z)+1)( C_{\nu_n}(z)+1) $ converges uniformly to $1$ on every compact of $\Delta_{\alpha,\beta_3}$. So, since $(\la +1)U(z)$ is equivalent to $z$ as $z$ tends to zero and since $\ds C_{\nu_n}=U\lf(z/H_{\nu_n}^{-1}(z)-1\ri)$, it suffices to show that   for all $z\in \Delta_{\alpha,\beta_3}$, \[\lim_{n\to\infty}\f{\f{z}{H_{ \nu_n}^{-1}(z)}-1
}{1-\f{z}{H_{ \nu_n}(z)}}=1.\]  
As in the proof of $(i)$, let us choose $\beta <\beta_3/\! 2$ \st for all $n$, $H_{\nu_n}^{-1}$ is defined on $\Delta_{\alpha',\beta}$,  $H_{\nu_n}(\Delta_{\alpha,\beta})\cup H_{\nu_n}^{-1}(\Delta_{\alpha,\beta})\subset \Delta_{\alpha',2\beta},$ and $$\lim_{n\to\infty}\f{H_{\nu_n}^{-1}(z)}{z}=1{\textrm{ uniformly on every compact of $\Delta_{\alpha',2\beta}$}}.$$ By analycity of the $H_{\nu_n}^{-1}$'s, the last assertion implies that$$\lim_{n\to\infty}(H_{\nu_n}^{-1})'(z)
=1{\textrm{ uniformly on every compact of $\Delta_{\alpha',2\beta}$}}.$$  We have, for $z\in\Delta_{\alpha,\beta}$,   \begin{eqnarray*}\f{\f{z}{H_{ \nu_n}^{-1}(z)}-1
}{1-\f{z}{H_{ \nu_n}(z)}}-1&=&\f{\f{zH_{ \nu_n}(z)}{H_{ \nu_n}^{-1}(z)}-H_{ \nu_n}(z)
}{H_{ \nu_n}(z)-z}-1\\ 
&=&\ff{H_{ \nu_n}(z)-z}\int_{[z,H_{ \nu_n}(z)]}f_n'(\xi)\ud \xi, \end{eqnarray*} where $f_n$ is still the function defined by $$f_n(\xi)=\f{H_{\nu_n}^{-1}(\xi)H_{\nu_n}(z)}{H_{\nu_n}^{-1}(z)}-\xi.$$ $f_n'$ tends to zero as $n$ tends to infinity, uniformly on every compact of $\Delta_{\alpha',2\beta}$, so $$\f{\f{z}{H_{ \nu_n}^{-1}(z)}-1
}{1-\f{z}{H_{ \nu_n}(z)}}-1$$ tends to zero when $n$ tends to infinity, and the result is proved.
\end{pr}

In the following, we shall refer to {\it weak convergence} for sequences of  positive finite measures on the real line: it is the convergence for which the test functions are the continuous bounded functions. 
\begin{Th}[L\'evy-Kinchine formula, part 1]\label{levyrect1} Let $\mu$ be a symmetric \pro measure, $(\nu_n)$ be a sequence of symmetric \pro measures and $k_n$ a sequence of integers tending to infinity \st $\nu_n^{\arc k_n}$ converges weakly to $\mu$. Then there exists a symmetric positive finite measure $G$ \st \begin{enumerate}\item the sequence of positive finite measures $\lf(k_n\f{t^2}{1+t^2}\ud \nu_n (t)\ri)$ converges weakly to $G$,\item the rectangular $R$-transform of $\mu$ has an analytic continuation to the complement of the real nonnegative half line and is given by the formula \begin{equation}\label{formuleLVKrect}C_\mu(z)=z\int_\R\f{1+t^2}{1-zt^2}\ud G(t).\end{equation}\end{enumerate}
 Moreover,  $G$ is symmetric and is the only positive finite measure $F$ \st \[C_\mu(z)=z\int_\R\f{1+t^2}{1-zt^2}\ud F(t).\] 
\end{Th}

\begin{pr}\begin{enumerate}\item The sequence $( \nu_n )$ converges weakly to $\delta_0$. Indeed, for every $n$, $C_{\nu_n^{\arc k_n}}=k_nC_{\nu_n}$, and by \teo \ref{tesprof?}, we have\begin{itemize} \item[(a)] $\underset{\substack{x\to 0\\ x <0}}{\lim}k_nC_{\nu_n}(x)=0$ uniformly in $n$,  \item[(b)] There exists $\beta> 0$ \st the sequence $(k_nC_{\nu_n})$ converges pointwise  on $(-\beta,0)$. \end{itemize}So \begin{itemize} \item[(a)]$\underset{\substack{x\to 0\\ x<0}}{\lim}C_{\nu_n}(x)=0$ uniformly in $n$,  \item[(b)] there exists $\beta> 0$ \st the sequence $(C_{\nu_n})$ converges pointwise to $0=C_{\delta_0}$  on $(-\beta,0)$. \end{itemize} 
\item The sequence of positive finite measures $\lf(k_n\f{t^2}{1+t^2}\ud \nu_n (t)\ri)$ is tight. Indeed, for $y>0$, \begin{eqnarray*}\int_{[-1/\!y,1/\!y]^c}k_n\f{t^2}{1+t^2}\ud \nu_n (t)& \leq & 2\int_{t\in\R}\f{1+t^2}{y^{-2}+t^2}\cdot k_n\f{t^2}{1+t^2}\ud \nu_n (t)\\
& = & -2k_n((i/\!y)G_{\nu_n}(i/\!y)-1).
 \end{eqnarray*}
We used the symmetry of $\nu_n$ in the second line. Let $v_n$ be as in the previous lemma. For $y>0$ small enough, 
\[\int_{[-1/\!y,1/\!y]^c}k_n\f{t^2}{1+t^2}\ud \nu_n (t) \leq  -\f{2k_nC_{\nu_n}(-y^2)}{1+v_n(-y^2)}
\leq  4 \lf| k_nC_{\nu_n}(-y^2)\ri|,  \]
which tends to zero uniformly in $n$ when $y$ tends to zero, by tightness of the sequence $\lf(\nu_n^{\arc k_n}\ri)$.
\item The sequence  of positive finite measures $\lf(k_n\f{t^2}{1+t^2}\ud \nu_n (t)\ri)$ is bounded. Indeed, choose $y\in (0,1)$ is \st $-y^2$ is  in the domain of the $v_n$'s of the previous lemma and  $\underset{n\to\infty}{\lim}k_nC_{\nu_n}(-y^2) = C_\mu(-y^2)$. Note that for all $t\in \R$, we have $\f{t^2}{1+t^2}<1<y^{-2}$, so $\f{y^{-2}+t^2}{1+t^2}<2y^{-2}$, hence $$\f{t^2}{1+t^2}<2y^{-2}\f{t^2}{y^{-2}+t^2}.$$ So we have, for each $n$,  \begin{eqnarray*}\int_{t\in\R}
k_n\f{t^2}{1+t^2}\ud \nu_n (t) & \leq & 2y^{-2}\int_{t\in\R}
k_n\f{t^2}{y^{-2}+t^2}\ud \nu_n (t)\\
& = & -2y^{-2}k_n((i/\!y)G_{\nu_n}(i/\!y)-1)\\
& = & -2y^{-2}k_n\f{C_{\nu_n}(-y^2)}{1+v_n(-y^2)}\\ & \leq & 4y^{-2}\lf|k_nC_{\nu_n}(-y^2)\ri|,\end{eqnarray*}which is  bounded uniformly in $n$.

\item Let us now recall a few facts about the {\textit{Poisson integral}} of positive measures on the real line which integrate $1/(1+t^2)$. 
If $M$ is such a measure,  for $y<0$ and $x\in \R$, let us define \[P_y(M)(x)= \int_{t\in \R}\f{y}{y^2+(x-t)^2}\ud M(t).\]Then $(x+iy)\mapsto P_y(M)(x)$ is harmonic and determines the measure $M$ (\cite{donog}, chapter II, \teo II). 
\\
Moreover, an easy computation shows that for each positive symmetric measure $M$ on the real line that integrates $1/(1+t^2)$,  the Poisson integral $P_y(M)(x)$ is the imaginary part of $\int_\R\f{\sqrt{z}}{t^2z-1}\ud M(t)$ (with $z\notin [0,+\infty)$, $x+iy=1/\!\sqrt{z}$, as it will be until the end of this proof). Indeed, since $M$ is symmetric, $$\ds\int_\R\f{\sqrt{z}}{t^2z-1}\ud M(t)=\int_\R\f{\sqrt{z}(t\sqrt{z}+1)}{t^2z-1}\ud M(t)=\int_\R\f{\sqrt{z}(t\sqrt{z}+1)}{\sqrt{z}(t\sqrt{z}+1)(t-\ff{\sqrt{z}})}\ud M(t)=$$ $$=\int_\R\f{\ud M(t)}{t-\ff{\sqrt{z}}}=\int_\R\f{(t-x)+iy}{(t-x)^2+y^2}\ud M(t).$$
\\
Now let us compute the Poisson integral of the measures $k_nt^2\ud \nu_n (t)$. Let $\alpha, \beta>0$ and $(v_n)$ be as in the previous lemma, 
$z\in \gab$. We have   \begin{eqnarray*} P_y(k_nt^2\ud \nu_n(t))(x) & = & \Im\lf(\int_\R\f{k_n\sqrt{z}t^2}{t^2z-1}\ud \nu_n(t)\ri).\end{eqnarray*}But since $\nu_n$ is symmetric, we have $$\int_\R\f{k_n\sqrt{z}t^2}{t^2z-1}\ud \nu_n(t)=k_n\int_\R\f{t(1+\sqrt{z}t)}{(t\sqrt{z}-1)(t\sqrt{z}+1)}\ud \nu_n(t)=k_n\int_\R\f{t\ud \nu_n(t)}{t\sqrt{z}-1},$$ which is equal, by an easy computation, to $$-k_n\f{(1/\!\sqrt{z})G_{\nu_n}(1/\!\sqrt{z})-1}{\sqrt{z}}.$$
So   $ \ds P_y(k_nt^2\ud \nu_n(t))(x)  =  -\Im\lf(k_n\f{C_{\nu_n}(z)}{\sqrt{z}(1+v_n(z))}\ri),$ which tends to   the imaginary part of $-C_\mu(z)/\!\sqrt{z}$, because $\nu_n^{\arc k_n}$ converges weakly to $\mu$ and $\lim\limits_{n\to\infty}v_n(z)=0$. 

The sequence  $\lf(k_n\f{t^2}{1+t^2}\ud \nu_n (t)\ri)$, bounded and tight, is relatively compact in the set of finite positive measures in the real line endowed with the topology of weak convergence (i.e. the topology defined by bounded continuous functions). 
If two measures $G,H$ are the weak limit of  subsequences of 
$\lf(k_n\f{t^2}{1+t^2}\ud \nu_n (t)\ri)$, then the measures 
$(1+t^2)\ud G(t)$ and $(1+t^2)\ud H(t)$  have the same Poisson integral on 
$\sqrt{\gab}$. Indeed, for $z\in\gab$, 
$$\ds P_y(k_nt^2\ud \nu_n(t))(x) =\int_\R\underbrace{\f{y(1+t^2)}{y^2+(x-t)^2}}_{\substack{{\textrm{continuous}}\\ {\textrm{bounded fct of $t$}}}}\f{k_nt^2}{1+t^2}\ud \nu_n(t)$$ 
tends at the same time to $P_y((1+t^2)\ud G(t))(x)$, to $P_y((1+t^2)\ud H(t))(x)$, and  to  the imaginary part of $-C_\mu(z)/\!\sqrt{z}$. It implies, by harmonicity, that they have the same Poisson integral on the lower half plane, which implies $H=G$. 
So the sequence $\lf(k_n\f{t^2}{1+t^2}\ud \nu_n (t)\ri)$ converges weakly to  a measure $G$, \st the Poisson integral $P_y((1+t^2)\ud G(t))(x)$,  is equal to  the imaginary part of $-C_\mu(z)/\!\sqrt{z}$.  Thus, the functions $$C_\mu(z)/\!\sqrt{z}\;\textrm{ and }\;\ds \int_\R \f{\sqrt{z}(t^2+1)}{1-t^2z}\ud G(t)$$ have the same imaginary part. For $z\in (-\infty,0)$, it follows that $$C_\mu(z)\;\textrm{ and }\;\ds z\int_\R \f{t^2+1}{1-t^2z}\ud G(t)$$  have the same real part, so, by analycity and since both tend to zero as $z$ goes to zero, they are equal. 
\item If $F$ is another positive finite measure \st 
$\ds C_\mu(z)=z\int_\R \f{t^2+1}{1-t^2z}\ud F(t)$, then 
$$\ds z\int_\R \f{t^2+1}{1-t^2z}\ud G(t)=z\int_\R \f{t^2+1}{1-t^2z}\ud F(t).$$ 
After division by $-\sqrt{z}$ and extraction of the imaginary part, this gives the equality of the Poisson integrals of 
$(1+t^2)\ud G(t)$ and of  $(1+t^2)\ud F(t)$, which implies $G=F$.\end{enumerate}
\end{pr}

The previous \teo implies that for all  $\arc$-infinitely divisible distribution $\mu$, there exists a unique positive finite measure $G$ \st $C_\mu$ is given by  equation (\ref{formuleLVKrect}). $G$ is symmetric (as limit of symmetric measures) and will be  called the {\textit{L\'evy measure}} of $\mu$. By injectivity of the rectangular $R$-transform, two different \pro measures cannot have the same L\'evy measure.
\begin{Th}[L\'evy-Kinchine formula, part 2]\label{levyrect2} Every symmetric positive finite measure on the real line is the L\'evy measure of a   $\arc$-infinitely divisible distribution. 
\end{Th}

Before the proof of the \teov let us state two lemmas.  The first one is about the rectangular $R$-transform of the symmetric Bernoulli distribution.
\begin{lem}\label{rect.R.tr.bernoulli}There exists a sequence $(\alpha_k)_{k\geq 2 }$ \st the associated power series has a positive radius of convergence and \st the rectangular $R$-transform with ratio $\la$ of $(\delta_1+\delta_{-1})/2$ is given by 
the formula $$C_{(\delta_1+\delta_{-1})/2}(z)=z+\sum_{k\geq 2}\alpha_k z^k.$$
\end{lem}

\begin{pr}By  the subsection called "The case of compactly supported \pro measures" of the section called "The rectangular $R$-transform" of \cite{fbg.AOP.rect} applied to $\mu=(\delta_1+\delta_{-1})/2$, we know that $$C_{(\delta_1+\delta_{-1})/2}(z)=c_2((\delta_1+\delta_{-1})/2)z+\sum_{k\geq 2}c_{2k}((\delta_1+\delta_{-1})/2) z^k,$$where the power series has a positive radius of convergence. So it suffices to prove that $c_2((\delta_1+\delta_{-1})/2)=1$, which follows from  the equation (\ref{PlayItSam}) of the present paper.
\end{pr}

We will also need a result about the way  dilation of \pro measures modify the rectangular $R$-transform. For $c>0$, let us denote by $D_c: x\mapsto cx$. For any distribution $\mu$, $D_c(\mu)$ is the push-forward of $\mu$ by $D_c$, i.e. $D_c(\mu): B\mapsto \mu(c^{-1}B)$. 

\begin{lem}\label{action.dilatation}For all $\mu$ symmetric \pro measure, for all $c>0$, \begin{equation}\label{27.11.03.1}C_{D_c(\mu)}(z)=C_\mu(c^2z). \end{equation}\end{lem}\begin{pr}
We have $\qquad G_{D_c(\mu)}
=
\ff{c}G_\mu(\f{z}{c}),$
\[{\textrm{so }}\qquad H_{D_c(\mu)}(z)=\f{\la}{c^2} \lf(G_\mu (\ff{c\sqrt{z}})\ri)^2+\f{(1-\la)c\sqrt{z}}{c^2}G_{\mu} (\ff{c\sqrt{z}})=\ff{c^2} H_\mu\lf(c^2z\ri),\] \begin{eqnarray}
{\textrm{i.e. }}\qquad\qquad H_{D_c(\mu)}&=&D_{\ff{c^2}}\circ H_\mu\circ D_{c^2},\nonumber\\ 
\qquad\qquad H_{D_c(\mu)}^{-1}&=&D_{\ff{c^2}}\circ H_\mu^{-1}\circ D_{c^2},\nonumber\\ 
{\textrm{then 
}}\qquad\qquad C_{D_c(\mu)}(z)&=&U\lf(\f{c^2z}{H^{-1}_\mu(c^2z)}-1\ri),\nonumber\\ {\textrm{that is }}\qquad\qquad C_{D_c(\mu)}(z)&=&C_\mu(c^2z). \nonumber\end{eqnarray} \end{pr}

\begin{prth} Let us denote by $\mc{M}$ the set of symmetric positive finite measures $G$ on the real line \st there exists a symmetric distribution $\mu$ whose rectangular $R$-transform  is given by  equation (\ref{formuleLVKrect}). We will show that $\mc{M}$ is the set of symmetric positive finite measures, proving  that $c\delta_0$ and $c(\delta_u+\delta_{-u})\in \mc{M}$ for all $c,u>0$, that $\mc{M}$ is stable under addition, and that $M$ is closed under weak convergence. Note that once this result is proved, it will be clear that any symmetric \pro measure with rectangular $R$-transform given by equation (\ref{formuleLVKrect}) will be  $\arc$-infinitely divisible. Indeed, denoting $$\begin{array}{cr}\ds C^{(G)}(z)=z\int_\R\f{1+t^2}{1-zt^2}\ud G(t)&\quad\quad 
(G\in \mc{M}),\end{array}$$ we have  $C^{(G)}=nC^{\lf(\f{G}{n}\ri)}.$
\begin{enumerate}\item For every $c>0$, $c\delta_0\in\mc{M}$. Indeed, by equation (\ref{27.11.03.1}), if $C^{(\delta_0)}=C_\mu$, then for every $c>0$, $C^{(c\delta_0)}=C_{\mu'}$, with $ \mu'=D_{c^{1/2}}(\mu)$, 
so it suffices to show that there exists a symmetric distribution whose rectangular $R$-transform is $C^{(\delta_0)}$. This distribution will appear as the limit in the {\textit{rectangular free central limit theorem}}: the sequence $D_{n^{-1/2}}\lf((\delta_1+\delta_{-1})/2)^{\arc n}\ri)$ converges weakly to a distribution with  rectangular $R$-transform  $C^{(\delta_0)}$  (we will see in the following that it stays true if one replaces $(\delta_1+\delta_{-1})/2$ by any symmetric \pro measure with variance equal to $1$). Indeed, let $C_n$ denote the rectangular $R$-transform of $D_{n^{-1/2}}\lf((\delta_1+\delta_{-1})/2)^{\arc n}\ri)$. By theorem \ref{tesprof?}, we have to prove that  \begin{itemize}\item[(a)] $\underset{\substack{x\to 0\\ x\in (-\infty,0)}}{\lim}C_{n}(x)=0$ uniformly in $n$,  \item[(b)] there exists $\beta> 0$ \st for all $y\in (0,\beta]$, the sequence $(C_{n}(-y))$ converges to $-y$. \end{itemize} Note that $C_n(z)=nC_{(\delta_1+\delta_{-1})/2}(z/n)$ (we used lemma \ref{action.dilatation} and the additivity of the rectangular $R$-transform (Theorem \ref{additivity.21.09.06})). Hence lemma \ref{rect.R.tr.bernoulli} allows to conclude.

\item For all $c,u>0$,  $c(\delta_u+\delta_{-u})\in \mc{M}$.  Indeed, we have \[C^{(c(\delta_u+\delta_{-u}))}(z)=2c\f{z(1+u^2)}{1-u^2z}=2\f{c(1+u^2)}{2u^2}\f{(u^2z)(1+1^2)}{1-(u^2z)}=C^{\lf(c'(\delta_1+\delta_{-1})\ri)}(u^2z),\] where $c'=\f{c(1+u^2)}{2u^2}$. So, by equation (\ref{27.11.03.1}), it suffices to show that for all $c>0$, there exists a distribution whose  rectangular $R$-transform is $C^{(c(\delta_1+\delta_{-1}))}$. It is the same to prove that there exists a distribution whose  rectangular $R$-transform is $C^{(\f{c}{4}(\delta_1+\delta_{-1}))}$.  This distribution will appear as the limit in the {\textit{rectangular free Poisson limit theorem}}: the sequence $\nu_n^{\arc n}$, with $\nu_n= \lf(1-\f{c}{n}\ri)\delta_0+\f{c}{2n}\lf(\delta_1+\delta_{-1}\ri)$, converges weakly to a distribution with  rectangular $R$-transform  $C^{(\f{c}{4}(\delta_1+\delta_{-1}))}$. 

Indeed, $G_{\nu_n}(z)= \f{z^2-1+c/n}{z(z^2-1)}$, so, if $(v_n)$ is a sequence of functions on $\gab$ as in the lemma \ref{27.11.03.3}, we have \[C_{\nu_n}(z)=\lf(\ff{\sqrt{z}}G_{\nu_n}\lf(\ff{\sqrt{z}}\ri)-1\ri)(1+v_n(z))=\f{cz}{n(1-z)}(1+v_n(z)),\]  so for $\mu_n=\underbrace{\nu_n\arc\cdots \arc\nu_n}_{\textrm{$n$ times}}$, \[C_{\mu_n}(z)=nC_{\nu_n}(z)=\f{cz}{1-z}(1+v_n(z)).\]
So by the properties of the functions $v_n$, we have both \[\lim_{\substack{z\to 0\\ |\arg z-\pi|<\alpha}}C_{\mu_n}(z)=0{\textrm{ uniformly in $n$}}\] and \[\forall z\in \gab, \lim_{n\to\infty}C_{\mu_n}(z)=\f{cz}{1-z}=C^{(\f{c}{4}(\delta_1+\delta_{-1}))}(z).\]
So, by  \teo \ref{tesprof?}, we know that there exists a distribution whose  rectangular $R$-transform is $C^{(\f{c}{4}(\delta_1+\delta_{-1}))}$.
\item $\mc{M}$ is stable under addition because $C_{\mu}+C_{\nu}=C_{\mu\arc \nu}$.
\item $\mc{M}$ is closed under weak convergence: let $(G_n)$ be a sequence of $\mc{M}$ that converges to a  finite measure $G$. Then  clearly, the sequence $\lf(C^{(G_n)}\ri)$ converges pointwise to $C^{(G)}$.  So, by \teo \ref{tesprof?}, to prove that $G\in\mc{M}$, it suffices to show that  \[\underset{\substack{x\to 0\\ x<0}}{\lim}C^{(G_n)}(x)=0 {\textrm{ uniformly in $n$.}}\]
For each $n$ and $x\in (0,1)$, since $G_n$ is symmetric,  $ C^{(G_n)}(-x^2)=-\int_{\R}\f{x^2+t^2x^2}{1+t^2x^2}\ud G_n(t)$, \[\forall t\in\R,  \f{x^2+t^2x^2}{1+t^2x^2}\leq\begin{cases}
x(x+1)& {\textrm{ if $-1/\!x^{1/\! 2}\leq t\leq 1/\!x^{1/\! 2}$,}}\\  1 & {\textrm{ otherwise.}}\end{cases}\] So \[\lf| C^{(G_n)}(x)\ri| \leq x(x+1)G_n(\R)+G_n\lf(\R-[-1/\!x^{1/\!2},1/\!x^{1/\! 2}]\ri), \] which tends to zero uniformly in $n$ when $x$ tends to $0$, by  boundedness and tightness of $\{G_n\ste n\in \n\}$.
\end{enumerate}
\end{prth}

Both previous theorems together allow us to state the following corollary. 
\begin{cor} A symmetric \pro measure $\mu$ is $\arc$-infinitely divisible \ssi there exists a sequence $(\nu_n)$ of symmetric \pro measures and a sequence $(k_n)$ of integers tending to infinity \st the sequence $\lf(\nu_n^{\arc k_n}\ri)$ tends to $\mu$. 
\end{cor}

\section{Rectangular Bercovici-Pata bijection}
In this section, we will show that the bijective correspondence between classical symmetric infinitely divisible distributions and rectangular  free infinitely divisible distributions is a homeomorphism, and that there exists a correspondence between limit theorems for sums of independent symmetric random variables and sums of free rectangular random variables. 

Let us recall a few facts about symmetric $*$-infinitely divisible distributions, that can be found in \cite{gne} (or \cite{feller2}, \cite{petrov} ... ). A symmetric \pro measure $\mu$ on the real line is  $*$-infinitely divisible \ssi there exists a finite positive symmetric measure $G$ \st          \[
\forall \xi\in \R, \int_{t\in\R}e^{it\xi}\ud \mu(t)=\exp\lf(\int_{t\in\R}(cos(t\xi)-1)\f{1+t^2}{t^2}\ud G(t)\ri).\]In this case, such a measure $G$ is unique, and we will call it the {\textit{L\'evy measure}} of $\mu$, and a sequence of symmetric \cids converges weakly \ssi the sequence of the corresponding L\'evy measures converges weakly. Moreover, in this case, the L\'evy measure of the limit will be the limit of the L\'evy measures.  

We can then define the {\textit{rectangular Bercovici-Pata bijection with ratio $\la$}}, denoted by $\Lambda_\la$,  from the set of  symmetric $*$-infinitely divisible distributions to the set of $\arc$-infinitely divisible distributions, that maps a $*$-infinitely divisible distribution to the $\arc$-infinitely divisible distribution with the same L\'evy measure. Let $\mu$, $\nu$ be two $*$-infinitely divisible distributions with L\'evy measures $G$, $H$. Then the L\'evy measures of  $\mu *\nu$ and of $\Lambda_\la (\mu)\arc\Lambda_\la(\nu)$ are both $G+H$, so we have $$\ds\Lambda_\la (\mu*\nu)=\Lambda_\la (\mu)\arc\Lambda_\la(\nu) .$$

\begin{Th}\label{7.12.03.1} The rectangular Bercovici-Pata bijection with ratio $\la$ is a homeomorphism, which means that a sequence of $\arc$-infinitely divisible distributions converges weakly \ssi the sequence of the corresponding  L\'evy measures converges weakly, and in this case, the L\'evy measure of the limit is  the limit of the L\'evy measures. 
\end{Th}
 
\begin{rmq} Note that, for $G$ symmetric positive finite measure, the function $C^{(G)}(z)$ can also be written, by symmetry, $$\ds C^{(G)}(z)=\int_\R\f{z+t\sqrt{z}}{1-t\sqrt{z}}\ud G(t).$$  
\end{rmq}

\begin{pr} Since the rectangular $R$-transform $C_\mu$ with ratio $1$ of a symmetric distribution $\mu$ is linked to its \trv $\vfi_\mu$ by the relation $C_\mu(z)=\sqrt{z}\vfi_\mu(1/\!\sqrt{z})$ (see paragraph 5 of \cite{defconv} for the construction of the Voiculescu transform, and use the fact that for symmetric distributions, the L\'evy measure is symmetric to obtain $C_\mu(z)=\sqrt{z}\vfi_\mu(1/\!\sqrt{z})$), the previous remark and \teo 5.10 of \cite{defconv} shows that the map $\Lambda_1$ is the restriction of the ``usual'' Bercovici-Pata bijection to the set of symmetric distributions. It has been proved in \cite{steen2} that the Bercovici-Pata bijection is a homeomorphism. So the \teo is proved in the case where $\la=1$. But  for every $*$-infinitely divisible distribution $\mu$, the formula of the rectangular $R$-transform with ratio $\la$ of $\Lambda_\la(\mu)$ does not depend on $\la$, so \teo \ref{tesprof?} allows us to claim that all $\La_\la$'s are homeomorphisms.  
\end{pr}

The next \teo furthers the analogy between the free rectangular convolution and the classical convolution of symmetric measures. As recalled in Theorem 3.3 of \cite{appenice}, it is proved in \cite{gne} that when $(\nu_n)$ is a sequence of symmetric \pro measures on the real line and $(k_n)$ is a sequence of integers tending to infinity, the sequence $\lf(\nu_n^{* k_n}\ri)$ converges weakly to a \cid  \ssi the sequence $\lf(\f{k_nt^2}{1+t^2}\ud \nu_n(t)\ri)$ of positive finite measures converges weakly to its L\'evy measure. By the \teo    \ref{levyrect1}, we know that it will be the case if  the sequence $\lf(\nu_n^{\arc k_n}\ri)$ converges weakly to the image of the $*$-infinitely divisible distribution by the rectangular Bercovici-Pata bijection. The following \teo states the converse implication. So we have, for all $*$-infinitely divisible distributions  $\mu$, \begin{equation}\label{11.12.03.1}{\textrm{$\lf(\nu_n^{* k_n}\ri)$ converges to $\mu$ }}\Longleftrightarrow{\textrm{ $\lf(\nu_n^{\arc k_n}\ri)$ converges to $\La_\la(\mu)$}}\end{equation}
\begin{Th}\label{eq.19.09.06.1}Let  $(\nu_n)$ be a sequence of symmetric \pro measures on the real line and $(k_n)$ be a sequence of integers tending to infinity. The sequence $\lf(\nu_n^{* k_n}\ri)$ converges weakly to an \cid \ssi the sequence  $\lf(\nu_n^{\arc k_n}\ri)$ converges weakly to its image  by the rectangular Bercovici-Pata bijection with ratio $\la$.
\end{Th}

\begin{pr} By what precedes, it suffices to prove that if the sequence $\lf(\f{k_nt^2}{1+t^2}\ud \nu_n(t)\ri)$ of positive finite measures converges weakly to a finite measure $G$, then the sequence  $\lf(\nu_n^{\arc n}\ri)$ converges weakly to the \rid with L\'evy measure $G$. Assume the sequence $\lf(\f{k_nt^2}{1+t^2}\ud \nu_n(t)\ri)$ of positive finite measures to converge weakly to a finite measure $G$.  \begin{enumerate}\item The sequence $(\nu_n)$ converges weakly to $\delta_0$:\\
Indeed, for all $\eps>0$, as the function $t\mapsto \f{t^2}{1+t^2}$ is increasing on $\R^+$, we have $$\nu_n\lf([-\eps,\eps]^c\ri) \leq \f{1+\eps^2}{\eps^2} \int_\R\f{t^2}{1+t^2}\ud \nu_n(t), $$ which tends to zero as $n$ tends to infinity, because the sequence $\lf(\f{k_nt^2}{1+t^2}\ud \nu_n(t)\ri)$ is bounded.
\item We have pointwise convergence of the rectangular $R$-transforms: \\ Let $\alpha,\beta$ and $(v_n)$ be as in the lemma \ref{27.11.03.3}. On $\gab$, we have \begin{equation}\label{7.12.03.2}C_{\nu_n^{\arc k_n}}(z)=k_nC_{ \nu_n}(z)=k_n\lf(\ff{\sqrt{z}}G_{ \nu_n}\lf(\ff{\sqrt{z}}\ri)-1\ri)(1+v_n(z)),\end{equation} but we have seen in the proof of theorem \ref{levyrect1} that \[k_n\lf(\ff{\sqrt{z}}G_{ \nu_n}\lf(\ff{\sqrt{z}}\ri)-1\ri)=z\int_\R\underbrace{\f{t^2+1}{1-t^2z}}_{\substack{{\textrm{continuous}}\\ {\textrm{bounded fct of $t$}}}}\f{k_nt^2}{1+t^2}\ud \nu_n(t),\] so, by pointwise convergence of the sequence $(v_n)$ to zero, the rectangular $R$-transform of $ \nu_n^{\arc k_n}$ converges pointwise to $z\mapsto z\int_\R\f{t^2+1}{1-t^2z}\ud G(t)$ on the set $\gab$. 
\item We have $\underset{\substack{y\to 0\\ y>0}}{\lim}C_{\nu_n^{\arc k_n}}(-y^2)=0$ uniformly in $n$: \\ By equation (\ref{7.12.03.2}) and $(i)$ of lemma \ref{27.11.03.3}, it suffices to prove that $$\underset{\substack{y\to 0\\ y>0}}{\lim}k_n((i/\!y)G_{ \nu_n}(i/\!y)-1)=0\quad {\textrm{  uniformly in $n$,}}$$ that is, since $\nu_n$ is symmetric,    $$\underset{\substack{y\to 0\\ y>0}}{\lim}\int_\R\f{y^2+t^2y^2}{1+t^2y^2}\f{k_nt^2}{1+t^2}\ud \nu_n(t)=0\quad {\textrm{  uniformly in $n$.}}$$ When $y<1$, $t\mapsto \f{y^2+t^2y^2}{1+t^2y^2}$ is $\leq 1$ and is increasing on $[0,\infty)$, so we have, for every $T>0$,   
$$\int_\R\f{y^2+t^2y^2}{1+t^2y^2}\f{k_nt^2}{1+t^2}\ud \nu_n(t)
\leq \int_{[-T,T]^c}\f{k_nt^2}{1+t^2}\ud \nu_n(t)+\f{y^2+T^2y^2}{1+T^2y^2}\int_{\R}\f{k_nt^2}{1+t^2}\ud \nu_n(t).$$
 Now fix $\eps >0$, choose $T>0$ \st for all $n$, 
$\int_{[-T,T]^c}\f{k_nt^2}{1+t^2}\ud \nu_n(t)\leq \eps$. For $y$ small enough, 
$\ds\f{y^2+T^2y^2}{1+T^2y^2}\sup_n\int_{\R}\f{k_nt^2}{1+t^2}\ud \nu_n(t)$ is less than $ \eps $, which closes the proof.
\end{enumerate}
\end{pr}

The following corollary could have been proved with the equation (\ref{27.11.03.1}), but the proof we give is shorter and does not use any computations.
\begin{cor}The rectangular Bercovici-Pata bijection commutes with the dilations $D_c$, $c>0$.\end{cor}
\begin{pr}Let $\mu$ be a \cidp Let, for each $n\geq 1$, $\nu_n$ be a symmetric distribution \st $\nu_n^{* n}=\mu$. We have \begin{eqnarray*}\La_\la\circ D_c(\mu)&=&\La_\la\circ D_c\lf(\nu_n^{* n}\ri)\\ &=&\La_\la\lf(D_c(\nu_n)^{* n}\ri) \\
&=&\lim_{n\to\infty}D_c(\nu_n)^{\arc n}. \end{eqnarray*} But from equation (\ref{27.11.03.1}) and additivity of the rectangular $R$-transform, we know that $$\forall n \geq 1, D_c(\nu_n)^{\arc n}=D_c\lf(\nu_n^{\arc n}\ri), $$ so, by continuity of $D_c$, $$\La_\la\circ D_c(\mu)=D_c\lf(\lim_{n\to\infty}\nu_n^{\arc n}\ri)$$ which is $D_c \circ \La_\la(\mu)$ by equivalence (\ref{11.12.03.1}). \end{pr}

Let us  define the {\textit{$\arc$-stable distributions}} to be the symmetric distributions whose orbit under the action of the group of the dilations is stable under $\arc$. The previous corollary allows us to give the following one. \begin{cor} The rectangular Bercovici-Pata bijection exchanges symmetric \csts and   $\arc$-stable distributions. Moreover, the index of any $*$-stable distribution $\mu$  (i.e. the unique $\alpha\in (0,2]$ \st for all $n\geq 1$, $\mu^{*n}=D_{n^{\ff{\alpha}}}(\mu)$) is preserved, i.e. one has $\Lambda_\la(\mu)^{\arc n}=D_{n^{\ff{\alpha}}}(\Lambda_\la(\mu))$.\end{cor}

The theorem \ref{eq.19.09.06.1} has another surprising consequence, which concerns classical \pro theory. It mights already be known by  specialists of limit theorems in classical \pro theory, but we since it can surprisingly be deduced from our results on Bercovici-Pata bijections, we state it and prove it here. In order to state it, we have to go further in the description of divisible distributions with respect to $\bxp$ and $*$: we have to give the L\'evy-Kinchine formulas for non symmetric infinitely divisible distributions. These distributions have been classified in \cite{defconv} and  \cite{gne}: a \pro measure on the real line $\mu$ is infinitely divisible with respect to $\bxp$ (resp. $*$) \ssi there exists a real number $\gamma$ and a positive finite measure on the real line $\sigma$ \st $\varphi_\mu(z)=\gamma+\int_\R\f{1+zt}{z-t}\ud \sigma(t)$ (resp. the Fourier transform is $\hat{\mu}(t)=\exp\lf[i\gamma t+\int_\R(e^{itx}-1-\f{itx}{x^2+1})\f{x^2+1}{x^2}\ud \sigma(x)\ri]$). Moreover, in this case, such a pair $(\gamma, \sigma)$ is unique, and we denote $\mu$ by $\nu_{\bxpp\!}^{\gamma, \sigma}$ (resp. $\nu_{\ast}^{\gamma, \sigma}$). Thus, one can define a bijection $\La$,  called the {\it Bercovici-Pata bijection}, from the set of \cids to the set of \fids by $$\La : \nu_{\ast}^{\gamma, \sigma}\mapsto \nu_{\bxpp\!}^{\gamma, \sigma}.$$ It is proved in \cite{appenice} that  for all  sequence $(\mu_n)$ of \pro measures and for all sequence $(k_n)$ of integers tending to $+\infty$, the sequence  $\mu_n^{\ast k_n}$ tends weakly to a \pro measure $\mu$ \ssi the sequence $\mu_n^{\bxp k_n}$ tends weakly to $\La(\mu)$. 
 By section \ref{20.09.06.1}, the infinitely divisible distributions with respect to $\bxp_1$ are the symmetric  infinitely divisible distributions with respect to $\bxp$ and the rectangular  Bercovici-Pata bijection with ratio $1$ is the restriction of the Bercovici-Pata bijection to the set of symmetric $*$-infinitely divisible distributions.
 
 \begin{cor} Let  $(\nu_n)$ be a sequence of symmetric \pro measures on the real line and $(k_n)$ be a sequence of integers tending to infinity. Let, for all $n$, $\rho_n$ be the push-forward of $\nu_n$ by the function $t\to t^2$. Then the sequence $\lf(\nu_n^{* k_n}\ri)$ converges weakly to a \pro measure \ssi the sequence  $\lf(\rho_n^{* k_n}\ri)$ converges weakly to a \pro measure. Moreover, this case, if one denotes  the L\'evy measure of the limit of $\lf(\nu_n^{* k_n}\ri)$ by $G$ (as the limit of such a sequence, the limit \pro measure has actually got to be symmetric and $*$-infinitely divisible), then  the limit of $\lf(\rho_n^{* k_n}\ri)$ is $\nu_{*}^{\gamma, \sigma}$, with $$\gamma=\int_{t\in \R}\f{1+t^2}{1+t^4}\ud G(t), \quad \sigma=\f{t^2+t}{t^2+1}\ud F(t),$$ where $F$ is the push-forward, by $t\to t^2$, of $G$. \end{cor}

\begin{pr}Let us first prove the equivalence. Recall that, as explained in section \ref{20.09.06.1},  for all $n$, the push-forward, by the function $t\to t^2$, of  $\nu_n^{\bxp_0 k_n}$ is  $\rho_n^{\bxp k_n}$. Hence for any symmetric \pro measure $\mu$,  if 
$$ \nu_n^{* k_n}\ninf \mu, \textrm{ i.e. } \nu_n^{\bxp_0 k_n}\ninf \La_0(\mu),$$ then 
$$ \rho_n^{\bxp k_n}\ninf \La_0(\mu)^2, $$
where $\La_0(\mu)^2$ denotes the push-forward, by the function $t\to t^2$, of $\La_0(\mu)$. Hence$$\rho_n^{* k_n}\ninf \La^{-1}(\La_0(\mu)^2).$$ Reciprocally, if there is a \pro measure $\rho$ on $[0+\infty)$ \st $$\rho_n^{* k_n}\ninf \rho, $$ then $$\rho_n^{\bxp k_n}\ninf \La(\rho), $$ i.e.  $\nu_n^{\bxp_0 k_n}$ converges weakly to the symmetric \pro measure $\nu$ whose push-forward by the square function is  $\La(\rho)$, which implies that  $\nu_n^{* k_n}$ converges weakly to a symmetric \pro measure. 

To prove the last part of the corollary, recall the fact from \cite{gne}, which is also recalled in \cite{appenice}, that for all  sequence $(\eta_n)$ of \pro measures on the real line, for all sequence $(k_n)$ of integers tending to $+\infty$, for all real number $a$ and all positive measure finite $H$ on the real line, we have the equivalence $$\eta_n^{* k_n}\ninf \nu_*^{a, H}\Longleftrightarrow \int_{t\in\R}\f{k_nt}{1+t^2}\ud \eta_n(t)\ninf a,\textrm{ and }
\f{k_nt^2}{1+t^2}\ud \eta_n(t)\ninf H,$$where he convergences of measures are with respect to the weak topology, i.e. against all continuous bounded functions (note that this equivalence could have been a way to prove the result without reference to the Bercovici-Pata bijections).  Suppose that $\lf(\nu_n^{* k_n}\ri)$ converges weakly to a \pro measure. This measure has to be symmetric and $*$-infinitely divisible. Let us denote its  L\'evy measure by $G$. Then   $$\int_{t\in\R}\f{k_nt}{1+t^2}\ud \rho_n(t)=\int_{t\in\R}\f{k_nt^2}{1+t^4}\ud \nu_n(t)=\int_{t\in\R}\f{1+t^2}{1+t^4}\f{k_nt^2}{1+t^2}\ud \nu_n(t)\ninf \int_{t\in\R}\f{1+t^2}{1+t^4}\ud G(t),$$ and for all continuous bounded function $f$,   $$\int_{t\in\R}f(t)\f{k_nt^2}{1+t^2}\ud \rho_n(t)=\int_{t\in\R}f(t^2)\f{k_nt^4}{1+t^4}\ud \nu_n(t)=\int_{t\in\R}f(t^2)\f{t^2+t^4}{1+t^4}\f{k_nt^2}{1+t^2}\ud \nu_n(t),$$which tends, when $n$ goes to infinity, to $$\int_{t\in\R}f(t^2)\f{t^2+t^4}{1+t^4}\ud G(t)=\int_{t\in\R}f(t)\f{t+t^2}{1+t^2}\ud F(t).$$ This concludes the proof.
\end{pr}

\section{Examples}\label{stzn.21.09.06}
In this section, we give examples of symmetric $*$-infinitely divisible distributions whose images by the rectangular Bercovici-Pata bijections we are able to give.  Unfortunately, there are as few examples as for the "classical"  Bercovici-Pata bijection. But in the section \ref{broderie+cafe=champion}, we shall give some matricial models for all $\arc$-infinitely divisible distributions.

\subsection{Rectangular Gaussian distribution and Marchenko-Pastur distribution}\label{6.2.05.1}
In  this section, we will identify the {\textit{rectangular Gaussian distribution}} $\nu$, that is the image, by the rectangular Bercovici-Pata bijection, of the Gaussian distribution with mean zero and variance one. The corresponding L\'evy measure is $\delta_0$, so the rectangular $R$-transform is $z$. We will show that unless $\la=0$, in which case $\nu=(\delta_{-1}+\delta_1)/2$,     $\nu$ is the symmetric distribution whose push forward by the function $x\to x^2$ has the density $$\f{\lf[4\la -(x-1-\la)^2\ri]^{1/2}}{2\pi \la x }\chi(x),$$where $\chi$ stands for the characteristic function of the interval $[(1-\la^{1/2})^2,(1+\la^{1/2})^2]$, which means that for all $n\geq 1$, the $2n$-th moment of $\nu$ is $1/\!\la$ times the $n$-th moment of the Marchenko-Pastur distribution with expectation $\la$ (the Marchenko-Pastur distributions 	are presented in section \ref{20.09.06.975}). 
 
Recall that the sequence $(c_{2n}(\mu))_{n\geq 1}$ of the {\it free cumulants with ratio $\la$} of a symmetric \pro measure $\mu$ with moments of any order, defined in the subsection called "The case of compactly supported \pro measures" of the section called "The rectangular $R$-transform" of \cite{fbg.AOP.rect}, are linked to the sequence $(m_{n}(\mu))_{n\geq 0}$ of its moments by the relation (see the proposition 3.5 of \cite{fbg.AOP.rect}): \begin{equation}\label{PlayItSam}\forall n\geq 1,\quad m_{2n}(\mu)=\ds\sum_{\pi\in \NCp(2n)}\la^{e(\pi)}\prod_{V\in \pi}c_{|V|}(\mu),\end{equation}where $\NCp(2n)$ is the set of noncrossing partitions of $\{1,\ldots, 2n\}$ in which all blocks have even cardinality, and where $e(\pi)$ denotes the number of blocks of $\pi$ with even minimum. 

The following lemma will be useful to study distributions coming from rectangular free \pro theory. A function $f$ defined on a conjugation-stable  subset of $\C$ is said to be {\it commuting with the conjugation} (abbreviated by  {\it c.w.c.}) if $f(\bar{z})=\overline{f(z)}$. Note that the function $z\to z^{1/2}$ is  c.w.c., whereas  $z\to \sqrt{z}$ is not.

\begin{lem}\label{29.04.04.1}If the   rectangular $R$-transform of a symmetric \pro measure $\mu$ extends to an analytic c.w.c. function in a neighborhood 
$B(0,r)$ of zero in the complex plane and tends to zero at zero, then the \pro measure has compact support, and the expansion of $C_\mu(z)$ for small $z$ is given by the formula\begin{equation}\label{9.12.03.1}C_\mu(z)=\sum_{n=1}^{+\infty}c_{2n}(\mu)z^n.\end{equation}
\end{lem}
\begin{pr} 
Let us define $T(z)=(\la z+1)(z+1)$. Note that $U$ is the inverse of $T-1$. Since the extension of $C_\mu$ tends to zero at zero, $z/\!H_\mu^{-1}(z)  $ extends to a neighborhood of zero \st we have, in this neighborhood, $$\f{z}{H_\mu^{-1}(z) }=T\lf(C_\mu(z)\ri),$$ and this function tends to $1$ at zero. Thus $H_\mu^{-1}(z)$ is one to one in a neighborhood of zero, and $H_\mu$ extends to an analytic  c.w.c. function in a neighborhood of zero \st $$\lim_{z\to 0}\f{H_\mu(z)}{z}=1.$$ 
So  the function $G_\mu(1/\!\sqrt{z})/\! \sqrt{z}$, which is equal to $$\f{\la-1+\lf[(1-\la)^2+4\la (H_\mu(z)/\! z)\ri]^{1/2}}{2\la}$$ if $\la>0$ and to $H_\mu(z)$ if $\la=0$, extends to an analytic c.w.c.  function in a neighborhood of zero. But since  $\mu$ is symmetric, for all $z$ in the complement of the real nonnegative half line, $$\f{G_\mu(1/\!\sqrt{z})}{\sqrt{z}}\!\!=\!\!\ff{\sqrt{z}}\int_\R\f{\ud \mu(t)}{\ff{\sqrt{z}}-t}\!\!
=\!\!\ff{2\sqrt{z}}\int_\R\lf[\ff{\ff{\sqrt{z}}-t}+\ff{\ff{\sqrt{z}}+t}\ri]\ud \mu(t)\!\!=\!\!\ff{z}\int_\R\f{\ud \mu(t)}{\ff{z}-t^2}\!\!=\!\!\ff{z}G_\rho(1/\!z),$$
where $\rho$ is the push forward of $\mu$ by the function $t\to t^2$. Hence the Cauchy transform of $\rho$ extends to an analytic c.w.c. function in a neighborhood of infinity. Thus, by the Stieltjes inversion  formula, $\rho$ is compactly supported, which implies that $\mu$ has compact support too.

$\mu$ has now been proved to be compactly supported. Then  the second part of the lemma, equation (\ref{9.12.03.1}), has been established in  the subsection called "The case of compactly supported \pro measures" of the section called "The rectangular $R$-transform" of \cite{fbg.AOP.rect}.
\end{pr}

So $\nu$ has compact support, and for all $n\geq 1$, $c_{2n}(\nu)=\delta_{1,n}$. 

Let us first treat the case where $\la=0$. By (\ref{PlayItSam}), all even moment of $\nu$ are $1$, so $\nu=(\delta_{-1}+\delta_1)/2$.

 Assume $\la >0$. By (\ref{PlayItSam}), the moments of $\nu$ are given by \begin{equation*}\forall n\geq 1,\;\; m_{2n}(\nu)=\sum_{\pi}\la^{e(\pi)}=\la^n\sum_{\pi}\lf(\ff{\la}\ri)^{o(\pi)},\end{equation*} where the sums are taken over noncrossing pairings of $\{1,\ldots,2n\}$ (a {\it noncrossing pairing} is a noncrossing partition where all classes have cardinality two, recall also that for a partition $\pi$, $e(\pi)$ and $o(\pi)$ are respectively the number of classes of $\pi$ with even and odd minimum). 
\begin{lem}Let $I=\{x_1<\cdots <x_n\}$ and $J=\{y_1<z_1<y_2<z_2<\cdots <y_n<z_n\}$ be  totally ordered sets. There is a bijection $\pi\to \ddot{\pi}$ from the set of noncrossing partitions of $I$ to the set of  noncrossing pairings of $J$ \st for all $\pi$, $$|\pi|=o(\ddot{\pi}).$$
\end{lem}

\begin{pr}Let us first construct the map $\pi\to \ddot{\pi}$ by induction on $n$, using the following well known result : a partition $\pi$ of a finite totally ordered set is noncrossing if and only if one of its classes $V$ is an interval and $\pi\backslash \{V\}$ is noncrossing (page 3 of \cite{spei98}). Consider a noncrossing partition $\pi$ of $I$. If $\pi$ has only one class, we define $\ddot{\pi}$ to be $$\{\{y_1,z_n\},\{z_1,y_2\},\{z_2,y_3\},\ldots,\{z_{n-1},y_n\}\}.$$ In the other case, a strict class $V$ of $\pi$ is an interval, $V=\{x_{k},x_{k+1}, \ldots, x_l\}$. Then we define   $\ddot{\pi}$ to be $$\ddot{\sigma}\cup \{\{y_{k},z_l\},\{z_k,y_{k+1}\},\{z_{k+1},y_{k+2}\},\ldots,\{z_{l-1},y_l\}\},$$ where $\ddot{\sigma}$ is the image (defined by the induction hypothesis) of the partition $$\sigma=\pi-\{V\}$$ of $I-V$  (it is easy to see that the result does not depend on the choice of the interval $V$).

The relation $|\pi|=o(\ddot{\pi})$ follows from the construction of $\pi\to \ddot{\pi}$.

Let us now prove, by induction on $n$, that $\pi\to \ddot{\pi}$ is a bijection. If  $n=1$, the result is obvious. Suppose the result to be proved to the ranks $1,\ldots, n-1$, and consider a noncrossing pairing  $\tau$ of $J$. Let us prove that there exists exactly one noncrossing partition $\pi$ of $I$ \st $  \ddot{\pi}=\tau$. Consider $l\in [n]$ minimal \st there exists $k<l$ \st $\{y_k,z_l\}$ is a class of $\tau$ (such an $l$ exists because it is the case of $n$). Then  it is easy to see that  $\{z_k,y_{k+1}\},\{z_{k+1},y_{k+2}\},\ldots,\{z_{l-1},y_l\}$ are classes of $\tau$, and any partition  $\pi$ of $I$ \st $  \ddot{\pi}=\tau$ must satisfy  $V:=\{x_{k},x_{k+1}, \ldots, x_l\}\in \pi$, and $$\ddot{\sigma}=\tau-\{\{y_{k},z_l\},\{z_k,y_{k+1}\},\{z_{k+1},y_{k+2}\},\ldots,\{z_{l-1},y_l\}\},$$ where $\sigma=\pi-\{V\}$ (partition of $I-V$). Thus, by the induction hypothesis, there exists  exactly one noncrossing partition $\pi$ of $I$ \st $  \ddot{\pi}=\tau$. 
\end{pr}

So 
the moments of $\nu$ are given by \begin{equation*}\forall n\geq 1,\;\; m_{2n}(\nu)=\la^n\ds\sum_{\pi\in \NC(n)}\lf(\ff{\la}\ri)^{|\pi|}.\end{equation*}
But for all $n\geq 1$, $\sum_{\pi\in \NC(n)}\lf(1/\!\la\ri)^{|\pi|}$ is the $n$-th moment of a distribution with all free cumulants being equal to $1/\!\la$, i.e. of the Marchenko-Pastur distribution with parameter $1/\! \la$ (see section \ref{20.09.06.975}). Thus  
the push-forward of $\nu$ by $t\to t^2$ is the push-forward of the Marchenko-Pastur distribution  with parameter $1/\! \la$ by the map $t\to \la t$, and has density  $$\f{\lf[4\la -(x-1-\la)^2\ri]^{1/2}}{2\pi \la x }\chi(x),$$ where $\chi$ stands for the characteristic function of the interval $[(1-\la^{1/2})^2,(1+\la^{1/2})^2]$. 
Hence we have proved the following result:\begin{Th}\label{identification.gaussienne.dujeudi}
The rectangular Gaussian distribution $\nu$ with ratio $\la$ has cumulants given by $$\forall n\geq 1, c_{2n}(\nu)=\delta_{n,1}.$$When $\la=0$, $\nu=(\delta_1+\delta_{-1})/2$. When $\la>0$, $\nu$ has density $$\f{\lf[4\la -(x^2-1-\la)^2\ri]^{1/2}}{2\pi \la |x| }\chi(x^2).$$ Its support is $[-1-\la^{1/2},-1+\la^{1/2}]\cup [1-\la^{1/2},1+\la^{1/2}] $.\end{Th}

Note that when $\la=1$, it is the well-known semi-circle law with radius two.

\begin{rmq}[Moments of the Marchenko-Pastur distribution] Note that the previous lemma, used with the fact that the free cumulants of the Marchenko-Pastur distribution with parameter $a$ are all  equal to $a$ (see \cite{hiai} p. 65),  gives us a formula for the  $n$-th moment of the  Marchenko-Pastur distribution with parameter $a$: it is equal to $\sum_{\pi}a^{o(\pi)}$, where the sum is taken over all noncrossing pairings of $[2n]$. This formula, proved using a random matrix approach, appeared already in an unpublished paper of  Ferenc Oravecz and D\'enes Petz.\end{rmq}

\begin{rmq}[Growth of the support in the related semigroup]\label{growth.support.21.09.06.normal} Let us define, for a fixed $\la\in [0,1]$ and for all $c>0$, $N_c$ to be the law of $c^\ff{2}X$, when $X$ is a random variable distributed according $\nu$ (i.e. $N_c$ is the image of the symmetric Gaussian law with variance $c$ by the rectangular Bercovici-Pata bijection). Then $(\{N_c\ste c>0\}, \arc)$ is an additive semigroup, whereas the size of the support of  $N_c$ is not linear in $c$ but in $c^\ff{2}$. This kind of phenomenon had already been observed in free probability.\end{rmq}

\subsection{Rectangular Cauchy distributions} \label{cauchysection.22.09.06}
This section could be called {\it missed appointment for the Cauchy distribution}. The Cauchy type, $\{ \mc{C}_t = \ff{\pi}\f{t\ud x}{x^2+t^2}\ste t>0\}$, is well known to be invariant under  many transformations. For example,  this set is the set of symmetric $*$- and $\bxp$-stable distributions with index $1$  ($\mc{C}_t$ has L\'evy measure $t\mc{C}_1$). But we are going to compute the set of $\arc$-stable distributions with index $1$, and it will appear that  unless $\la=1$, it is not the Cauchy type. 

So  let us fix $\la\in [0,1]$, and let us denote the image $\Lambda_\la(\mc{C}_t)$ of the symmetric Cauchy law $\mc{C}_t$ with index $t$ by the Bercovici-Pata bijection with ratio $\la$ by $\nu_t$. We know that the rectangular $R$-transform with ratio $\la$ of $\nu_t$ admits an analytic extension to $\C\backslash \R^+$ given by the formula $it\sqrt{z}$. So, by the remark \ref{21.09.06.1},  $H^{-1}_{\nu_t}$ admits an analytic extension to $\C\backslash (\R^+\cup\{-\ff{t^2},\ff{t^2\la^2}\})$ given by the formula  
$$H^{-1}_{\nu_t}(z)=\f{z}{T(C_{\nu_t}(z))}=\f{z}{(\la it\sqrt{z}+1)(it\sqrt{z}+1)}, $$ where $T(X)=(\la X+1)(X+1)$.

In order to compute $H_{\nu_t}(x)$, we have to invert the previous formula and to remember that  it is a bijection form a neighborhood of $0$ in $\C\backslash \R^+$ to a neighborhood of zero in $\C\backslash \R^+$, equivalent to  $x$ in the neighborhood of zero. So let us fix $x\in \C\backslash \R^+$, and denote $H_{\nu_t}(x)$ by $z$. If $x$ is closed enough from zero, we have $x=H^{-1}_{\nu_t}(z)$, hence  
\begin{eqnarray*}x(\la it\sqrt{z}+1)(it\sqrt{z}+1)&=&z\\ -(\la xt^2+1)z+itx(\la+1)\sqrt{z}+x&=&0,\end{eqnarray*}Note that $(itx(\la+1))^2-4x(-(\la xt^2+1))=4x-t^2x^2(\la-1)^2$, so
$$\sqrt{z}=\f{itx(\la+1)\pm \sqrt{4x-t^2x^2(\la-1)^2}}{2\la xt^2+2}.$$ But when $x$ goes to zero in  $\C\backslash \R^+$, $$\sqrt{4x-t^2x^2(\la-1)^2}\sim 2\sqrt{x},\textrm{ so }\f{itx(\la+1)\pm \sqrt{4x-t^2x^2(\la-1)^2}}{2\la xt^2+2}\sim \pm \sqrt{x},$$ hence 
$$\sqrt{z}=\f{itx(\la+1)+ \sqrt{4x-t^2x^2(\la-1)^2}}{2\la xt^2+2}=\f{itx(\la+1)+ \sqrt{x}(4-t^2x(\la-1)^2)^\ff{2}}{2\la xt^2+2},$$ hence $$H_{\nu_t}(x)=\left( \f{itx(\la+1)+ \sqrt{x}(4-t^2x(\la-1)^2)^\ff{2}}{2\la xt^2+2}\right)^2=\f{x}{4}\left( \f{it\sqrt{x}(\la+1)+ (4-t^2x(\la-1)^2)^\ff{2}}{\la xt^2+1}\right)^2.$$ 
But by the remark \ref{21.09.06.1}, $$\ff{\sqrt{x}}G_{\nu_t}(\ff{\sqrt{x}})=V(\f{H_{\nu_t}(x)}{x})=\f{\la-1+[(\la-1)^2+4\la H_{\nu_t}(x)/x]^\ff{2}}{2\la}.$$
So we compute 
\begin{eqnarray*}&(\la-1)^2+4\la H_{\nu_t}(x)/x&\\ &=(\la-1)^2+ \f{\la}{(\la xt^2+1)^2}\{it\sqrt{x}(\la+1)+ [4-t^2x(\la-1)^2]^\ff{2}\}^2&\\ 
&=(\la-1)^2+ \f{\la}{(\la xt^2+1)^2}\{-t^2x(\la+1)^2+4-t^2x(\la-1)^2+2it\sqrt{x}(\la+1)[4-t^2x(\la-1)^2]^\ff{2}\}& \end{eqnarray*}
Hence
\begin{eqnarray*}& (\la xt^2+1)^2[(\la-1)^2+4\la H_{\nu_t}(x)/x]&\\
&=(\la-1)^2+2\la x t^2(\la-1)^2+\la^2x^2t^4(\la-1)^2-2\la xt^2(\la^2+1)+4\la+2it\la\sqrt{x}(\la+1)[4-t^2x(\la-1)^2]^\ff{2} & \\
&=(\la-1)^2+4\la+2\la xt^2[(\la-1)^2-\la^2-1]+\la^2x^2t^4(\la-1)^2+2it\la\sqrt{x}(\la+1)[4-t^2x(\la-1)^2]^\ff{2}  &\\
&=(\la+1)^2-4\la^2xt^2+\la^2x^2t^4(\la-1)^2+2it\la\sqrt{x}(\la+1)[4-t^2x(\la-1)^2]^\ff{2}  &\\
&=\{(\la+1)+it\la\sqrt{x}[4-t^2x(\la-1)^2]^\ff{2} \}^2&
\end{eqnarray*}

Hence for all $x\in \C\backslash \R^+$ closed enough to zero, \begin{eqnarray*}G_{\nu_t}(\ff{\sqrt{x}})&=&\f{\sqrt{x}}{2\la}\{\la-1+\underbrace{[(\la-1)^2+4\la H_{\nu_t}(x)/x]^\ff{2}}_{\f{(\la+1)+it\la\sqrt{x}[4-t^2x(\la-1)^2]^\ff{2}}{\la xt^2+1}}\}\\ 
&=&\f{\sqrt{x}}{2\la}\f{(\la-1)(\la xt^2+1)+\la+1+it\la\sqrt{x}[4-t^2x(\la-1)^2]^\ff{2}}{\la xt^2+1}\\ 
&=&\sqrt{x}\lf\{\f{2+(\la-1) xt^2+2it\sqrt{x}[1-t^2x(\la-1)^2/4]^\ff{2}}{2(\la xt^2+1)}\ri\}.
 \end{eqnarray*} 
So for all $z\in \C^-$ faraway enough from zero (in a non tangential way), for $x\in\C\backslash\R^+$ \st $z=1/\sqrt{x}$, i.e. $x=1/z^2$, $$G_{\nu_t}(z)=G_{\nu_t}\lf(\ff{\sqrt{x}}\ri)=\sqrt{x}\lf\{\f{2+(\la-1) xt^2+2it\sqrt{x}[1-t^2x(\la-1)^2/4]^\ff{2}}{2(\la xt^2+1)}\ri\}$$ $$=\ff{ z} \ff{\f{ 2\la t^2}{z^2}+2}\lf\{2+(\la-1)\f{ t^2}{z^2}+\f{2it}{z}\lf[1-\f{t^2(\la-1)^2}{4z^2}\ri]^\ff{2}\ri\}$$ $$=\f{z}{ \la t^2+z^2}\lf\{1+(\la-1)\f{ t^2}{2z^2}+\f{it}{z}\lf[1-\f{t^2(\la-1)^2}{4z^2}\ri]^\ff{2}\ri\}$$
Note that for $\rho \in \R^+$, for all $z\in\C\backslash\{0\}$,  $$1-\f{t^2(\la-1)^2}{4z^2}=-\rho\Leftrightarrow z^2=\f{t^2(1-\la)^2}{4(1+\rho)}\Leftrightarrow z=\pm\f{t(1-\la)}{2(1+\rho)^\ff{2}}\Rightarrow z\in[-\f{t(1-\la)}{2}, \f{t(1-\la)}{2}].$$So the function $z\mapsto \lf[1-\f{t^2(\la-1)^2}{4z^2}\ri]^\ff{2}$ extends analytically to $\C\backslash [-\f{t(1-\la)}{2}, \f{t(1-\la)}{2}]$ with the same formula.   So, by analycity, we have \begin{equation}\label{audition.20.04.1}\forall z\in\C^-\backslash \{i\la^\ff{2}t\}, \quad\quad G_{\nu_t}(z)=\f{z}{ \la t^2+z^2}\lf\{1+(\la-1)\f{ t^2}{2z^2}+\f{it}{z}\lf[1-\f{t^2(\la-1)^2}{4z^2}\ri]^\ff{2}\ri\}.\end{equation}

\begin{rmq}Note that $G_{\nu_t}(z)$ has to be analytic at $- i\la^{\ff{2}}t$.  If the formula we give had no analytic extension at $- i\la^{\ff{2}}t$, we would have made a mistake. Hopefully, one can check that the pole of $z\mapsto \f{z}{ \la t^2+z^2}$ at $-i\la^{\ff{2}}t$ is simple, and  that the function  $$z\mapsto 1+(\la-1)\f{ t^2}{2z^2}+\f{it}{z}\lf[1-\f{t^2(\la-1)^2}{4z^2}\ri]^\ff{2}$$ has a zero at  $- i\la^{\ff{2}}t$. \end{rmq}

 In order to compute $\nu_t$, we are going to use  the following lemma.

\begin{lem}\label{19.09.06.8}Let $\nu$ be  \pro measure on the real line \st the restriction of $G_\nu$ to $\C^-$ extends analytically to  an open set containing $\C^-\cup I$, where $I$ is an open interval. Then the restriction of $\mu$ to $I$ admits an analytic density: $x\in I \mapsto \ff{\pi}\Im G_\nu(x).$
\end{lem}

\begin{pr}Let us define, for $t\geq 0$, $$\rho_t : x\in \R\mapsto\begin{cases} \ff{\pi}\Im G_\nu(x-it)&\textrm{if $x\in I$,}\\ 0& \textrm{in the other case.}\end{cases}$$ Then for all $t>0$,
$\rho_t$ is well known to be the restriction, to $I$, of the density of $\nu*\mc{C}_t$. Moreover, $\nu*\mc{C}_t$ converges weakly (i.e. against any continuous bounded function)
to $\nu$ as $t$ tends to zero. So it suffices to prove that for all $f$ compactly supported continuous function on $I$,
$\int f(x)\rho_t(x)\ud x$ tends to $\int f(x)\rho_0(x)\ud x$ when $t$ goes to zero, which is an easy application of the dominated convergence theorem.\end{pr}

This lemma allows us to claim that the restriction of $\nu_t$ to $\R\backslash [-\f{t(1-\la)}{2}, \f{t(1-\la)}{2}]$ has an analytic density given by the function \begin{equation}\label{audition.21.04.1}x\mapsto \f{t}{ \pi(\la t^2+x^2)}\lf[1-\f{t^2(\la-1)^2}{4x^2}\ri]^\ff{2}\quad\quad(x\in \R\backslash [-\f{t(1-\la)}{2}, \f{t(1-\la)}{2}]).\end{equation}

In order to prove that $\nu_t$ is carried by  $\R\backslash [-\f{t(1-\la)}{2}, \f{t(1-\la)}{2}]$ and has this density, 
it suffices to prove that there is no mass out of  $\R\backslash [-\f{t(1-\la)}{2}, \f{t(1-\la)}{2}]$, i.e. that  $I=1$, with $$I:=\int_{\R\backslash [-\f{t(1-\la)}{2}, \f{t(1-\la)}{2}]}\f{t}{ \pi(\la t^2+x^2)}\lf[1-\f{t^2(\la-1)^2}{4x^2}\ri]^\ff{2}\ud x.$$We have \begin{eqnarray*}I&=&\f{2t}{\pi}\int_{\f{t(1-\la)}{2}}^{+\infty}\ff{ \la t^2+x^2}\lf[1-\f{t^2(\la-1)^2}{4x^2}\ri]^\ff{2}\ud x\\  
&=&\f{-2t}{\pi}\int_{\f{t(1-\la)}{2}}^{+\infty}\f{x^2}{ \la t^2+x^2}\lf[1-\f{t^2(\la-1)^2}{4x^2}\ri]^\ff{2}\ud \lf(\ff{x}\ri)\\
&=&\f{-2t}{\pi}\int_{\f{t(1-\la)}{2}}^{+\infty}\f{1}{ \f{\la t^2}{x^2}+1}\lf[1-\f{t^2(\la-1)^2}{4x^2}\ri]^\ff{2}\ud \lf(\ff{x}\ri)\\
&=&\f{2t}{\pi}\int_0^{\f{2}{t(1-\la)}}\f{1}{ \la t^2u^2+1}\lf[1-\f{t^2(\la-1)^2u^2}{4}\ri]^\ff{2}\ud u\\
&=&\f{2t}{\pi}\f{2}{(1-\la)t}\int_0^{\f{2}{t(1-\la)}}\f{1}{ \f{4\la}{(1-\la)^2} \lf(\f{(1-\la)tu}{2}\ri)^2+1}\lf[1-\f{t^2(\la-1)^2u^2}{4}\ri]^\ff{2}\ud \lf(\f{(1-\la)tu}{2}\ri)\\
&=&\f{4(1-\la)}{\pi}\int_0^{1}\f{1}{ 4\la v^2+(1-\la)^2}\lf[1-v^2\ri]^\ff{2}\ud v\end{eqnarray*}
Let us define $y=\f{v}{(1-v^2)^\ff{2}}$. When $v$ goes from $0$ to $1$, $y$ goes increasingly from $0$ to $+\infty$. One has $\ud y= \f{\ud v}{(1-v^2)^\f{3}{2}}$, and $v^2=\f{y^2}{1+y^2}$, hence $$\lf[1-v^2\ri]^\ff{2}\ud v= \lf[1-v^2\ri]^2\ud y=\f{\ud y}{(1+y^2)^2}$$ and $$\f{1}{ 4\la v^2+(1-\la)^2}=\f{1}{ 4\la \f{y^2}{1+y^2}+(1-\la)^2}=\f{1+y^2}{ 4\la y^2+(1-\la)^2(1+y^2)}=\f{1+y^2}{ (1+\la)^2y^2+(1-\la)^2}.$$So 
\begin{eqnarray*}I&=&2(1-\la)\int_\R\f{1}{ (1+\la)^2y^2+(1-\la)^2}\f{\ud y}{\pi(1+y^2)}\\
&=&\f{2(1-\la)}{ (1+\la)^2}\int_\R g(y)\f{\ud y}{\pi(1+y^2)},
\end{eqnarray*}
where $g(y)=\f{1}{ y^2+a^2}$, with $a=\f{1-\la}{1+\la}$. Note that $$g(y)=\ff{(y+ia)(y-ia)}=\ff{2ia}\lf(\ff{y-ia}-\ff{y+ia}\ri),$$ so the well known formula of the Cauchy transform of the Cauchy distribution with parameter $1$ gives us 
$$I=
\f{2(1-\la) }{ (1+\la)^2 } \cdot\ff{2ia} \lf( \ff{-i-ia}-\ff{i+ia}\ri) = \f{2(1-\la)}{ (1+\la)^2}\cdot\ff{2ia}\cdot\f{2i}{a+1}$$ $$=\f{2(1-\la)}{ (1+\la)^2}\cdot\f{1+\la}{1-\la}\cdot\ff{a+1}=\f{2}{(1+\la)a+(1+\la)}=1.$$
So we have proved the following result : 
\begin{propo}\label{domage.pr.cauchy.rect}
For all $\la\in [0,1]$, for all $t>0$, the image of the symmetric Cauchy law with parameter $t$ by the  Bercovici-Pata bijection with ratio $\la$ is $$1_{\R\backslash [-\f{t(1-\la)}{2}, \f{t(1-\la)}{2}]}(x) \f{t}{ \pi(\la t^2+x^2)}\lf[1-\f{t^2(\la-1)^2}{4x^2}\ri]^\ff{2}\ud x.$$Its support is $\R\backslash \lf(-\f{t(1-\la)}{2}, \f{t(1-\la)}{2}\ri)$
\end{propo}

The Cauchy type is well known to be invariant by the push-forward by the function $t\to 1/t$. In the following corollary, we are going to see that again, unless $\la=1$, things append differently for $\arc$-stable laws with index $1$.
\begin{cor}[Push-forward by the function $x\mapsto 1/x$]
\begin{itemize}\item[-] For $\la=1$,  the push-forward, by the function $x\mapsto 1/x$, of the measure presented in the proposition \ref{domage.pr.cauchy.rect} is the symmetric Cauchy law with parameter $1/t$. 
\item[-] For $\la\in [0,1)$, it is the  measure carried by $\lf[-\f{2}{t(1-\la)},\f{2}{t(1-\la)}\ri]$ with density \begin{equation}\label{compile.reussie.1}x\mapsto \ff{\pi(\la t^2x^2+1)}\lf[1-\f{x^2t^2(\la-1)^2}{4}\ri]^\ff{2}\quad\quad(x\in \lf[-\f{2}{t(1-\la)},\f{2}{t(1-\la)}\ri]).\end{equation} (a) When $\la=0$, it  is equal to \begin{equation}\label{compile.reussie.2}x\mapsto \f{t}{2\pi}\lf[\f{4}{t^2}-x^2\ri]^\ff{2}\quad\quad(x\in \lf[-\f{2}{t},-\f{2}{t}\ri]), \end{equation}hence it is  the symmetric semi-circle law with radius $2/t$ and variance $1/t^2$. \\
(b) When $\la\in (0,1)$, it is equal to  \begin{equation}\label{compile.reussie.3}x\mapsto \f{1-\la}{2\pi t\la}\ff{\ff{\la t^2}+x^2}\lf[\lf(\f{2}{t(1-\la)}\ri)^2-x^2\ri]^\ff{2}\quad\quad(x\in \lf[-\f{2}{t(1-\la)},\f{2}{t(1-\la)}\ri]).\end{equation}This density is the one of a symmetric semi-circle law with radius $\f{2}{t(1-\la)}$ and variance $\ff{t^2(1-\la)^2}$ times the density of a Cauchy law with parameter $\ff{t\la^\ff{2}}$ times $\f{\pi}{(1-\la)t^2\la^\ff{2}}$.\end{itemize}
\end{cor}

\begin{pr}Let $\rho$ be the density of a \pro on a Borel set $I$. Then for all bounded Borel function $f$,  $$\int_If\lf(\ff{x}\ri)\rho(x)\ud x=-\int_If\lf(\ff{x}\ri)\rho\lf(\ff{1/x}\ri)\ff{(1/x)^2}\ud \lf(\ff{x}\ri)=\int_Jf(y)\f{\rho(1/y)}{y^2}\ud y,$$where $J=\{1/x\ste x\in I\}$. Hence the push-forward of $\rho(x)\ud x$ by $x\mapsto 1/x$ is carried by $J$ and has density $\f{\rho(1/y)}{y^2}$. This proves the (well known) result for $\la=1$, and (\ref{compile.reussie.1}).  (\ref{compile.reussie.2}) and (\ref{compile.reussie.3}) follow easily. To recognize the products of semi-circle and Cauchy densities, just remember that for all $r>0$, the semi-circle law with radius $r$ has variance $r^2/4$ and density $$x\in [-r,r]\mapsto \f{2}{\pi r^2}\lf[r^2-x^2\ri]^\ff{2},$$ and that the Cauchy law with parameter $c>0$ is $ \mc{C}_c = \ff{\pi}\f{c\ud x}{x^2+c^2}$. 
\end{pr}

\subsubsection*{Question} Inspired by what happens in the cases $\la=0$ and $\la=1$, we ask the following question, the answer of which could have spared us the long previous calculus: is there, for each $\la$ in $[0,1]$, a functional $f_\la$ from the set of symmetric \pro measures on $\R$ into the set of \pro measures on $\R$ \st for all $\mu,\nu$ symmetric \pro measures, $\mu\arc\nu$ is the only symmetric \pro measure satisfying $$f_\la(\mu\arc\nu)=f_\la(\mu)\bxp f_\la(\nu)\;\; \textrm{ ?}$$ 
Note that in the case $\la=1$, the functional $f_\la(\mu)=\mu$ works, and in the case $\la=0$, the functional which maps a measure to its push-forward by the square function works. 

\subsection{Rectangular analogues of symmetrized Poisson distributions}\label{bon.vieux.dede.20.09.06}
\subsubsection{The general case $\la\in [0,1]$}\label{mille.petites.mains}Let us define the {\it symmetric Poisson distribution with parameter} $c>0$ to be the $*$-infinitely divisible distribution with L\'evy measure   $\f{c}{4}(\delta_1+\delta_{-1})$. It can also be seen as the law of $X-Y$, where $X,Y$ are independent random variables with (unsymmetric) Poisson law with parameter $c/\! 2$, or as the weak limit of   $$\lf(\lf(1-\f{c}{n}\ri)\delta_0+\f{c}{2n}(\delta_{-1}+\delta_1)\ri)^{* n}.$$  The  rectangular $R$-transform of  its image $P_c$ by the rectangular Bercovici-Pata bijection with ratio $\la$ is $$\f{cz}{1-z}.$$
Hence $$H_{P_c}^{-1}(z)=\f{z}{T(C_{P_c}(z))}=\f{z(1-z)^2}{[(\la c-1)z+1][(c-1)z+1]},$$ whose inversion would be very heavy because of a third degree equation. So we know the rectangular $R$-transform of this law, but we don't give the law (except when $\la=0$, see below). 

However, we know that for all $c,c'>0$, $P_c\arc P_{c'}=P_{c+c'}$, and in the section \ref{france.q=chinois.a.Belleville}, we shall give a nice matricial model for this distributions.

\subsubsection{The particular case $\la=0$}\label{20.09.06.975} Let us recall the definition, for $c>0$, of the  {\it Marchenko-Pastur law} with parameter $c$, also called {\it free Poisson law} with parameter $c$. It is the \pro measure on $[0,+\infty)$  $$M_c=\begin{cases}\f{(4c -(x-1-c)^2)^\ff{2}}{2\pi x}\chi_c(x)\ud x&\textrm{if $c\geq 1$,}\\ (1-c)\delta_0+\f{(4c -(x-1-c)^2)^\ff{2}}{2\pi x}\chi_c(x)\ud x&\textrm{if $0< c <1$,}\end{cases}$$ where $\chi_c$ is the characteristic function of  $[(1-c^\ff{2})^2,(1+c^\ff{2})^2]$. It is well known (see \cite{hiai}) that its Voiculescu transform is $cz/(z-1)$. 

Let us consider  the symmetric law whose push-forward by $t\to t^2$ is $M_c$. By section  \ref{20.09.06.1}, its rectangular $R$-transform with ratio $0$ is given b the formula $cz/(1-z)$. Hence 
for $\la=0$, the distribution $P_c$ introduced in the previous section \ref{mille.petites.mains}  is this distribution, hence (still for $\la=0$), $$P_c=\begin{cases}\f{(4c -(x^2-1-c)^2)^\ff{2}}{2\pi }\chi_c(x^2)\ud x&\textrm{if $c\geq 1$,}\\ (1-c)\delta_0+\f{(4c -(x^2-1-c)^2)^\ff{2}}{2\pi}\chi_c(x^2)\ud x&\textrm{if $0< c <1$.}\end{cases}$$

\begin{rmq}[Growth of the support in the semigroup] Note that we observe the same kind of phenomenon as in the remark \ref{growth.support.21.09.06.normal}: in the additive semigroup $(\{P_c\ste c>0\}, \arco\})$, the size of the support of $P_c$ is not linear in $c$, but  in $c^\ff{2}$.\end{rmq}

\section{A matricial model for the rectangular Bercovici-Pata bijection}\label{broderie+cafe=champion}
In the previous sections, the proofs rely on integral transforms and complex analysis. We will construct, in this subsection, a matricial model  for the $\arc$-infinitely divisible laws and present in a maybe more palpable way the Bercovici-Pata bijection with ratio $\la$. 

In this section, $d,d'$ will represent dimensions of rectangular matrices, because $n$ will be used to another role.
 For any distribution $\mathbb{P}$ and any function $f$ on a set of matrices, $\E_{\mathbb{P}}(f(M))$ denotes $\int f(M)\ud \mathbb{P}(M)$. Let us recall that the {\it singular law} of a matrix $M$ designates the uniform distribution on the spectrum of $|M|:=(MM^*)^{\ff{2}}$. Let us define the {\it symmetrization} $\tilde{\mu}$ of a distribution $\mu$ on the real line: it is the distribution which maps a Borel set $B$ to $(\mu(B)+\mu(-B))/2$. The symmetrization of the singular law of a matrix $M$ will be denoted by $\tilde{\mu}_{|M|}$.

We are going to construct, in the same way as in \cite{fbg.ID} and in \cite{cabduv.ID}, for each $d,d'\geq 1$, for each symmetric \cid $\mu$, an infinitely divisible  distribution \pmu on the set of \dtpd complex matrices \st for all $\mu,\nu$, $\mathbb{P}_{d,d'}^\mu*\mathbb{P}_{d,d'}^\nu=\mathbb{P}_{d,d'}^{\mu *\nu}$ and \st the symmetrization of the singular law  of  $M$ (with $M$ random matrix distributed according to $\mathbb{P}_{d,d'}^\mu$)  goes from $\mu$ to its image by the rectangular Bercovici-Pata bijection with ratio $\la$ when $d,d'\to\infty , \f{d}{d'}\to\la$.

Let us introduce the heuristic argument that led us to choose the model we will present.
Consider a symmetric \cid $\mu$, and two sequences $(\nu_n)$ (symmetric \pro measures), $(k_n)$ (integers tending to infinity) \st  $\ds\nu_n^{*k_n}$ tends weakly to $\mu$. 
Define, for all $1\leq d\leq d'$ and each $n\geq 1$,  $\mathbb{Q}_{d,d'}^{\nu_n}$ to be the law of the \dtpd random matrix $$\ds U\lf[X_{n,i}\delta_i^j\ri]_{\substack{1\leq i\leq d\\ 1\leq j\leq d'}}V$$ where $U$ (resp. $V$) is a uniform $d\tii d$ (resp. $d'\tii d'$) unitary random matrix, $X_{n,1},\ldots, X_{n,d}$ are distributed according to $\nu_n$, and $U,V, X_{n,1}\ldots$ are independent. 

Then if one fixes $n$ and lets $d,d'$ go to infinity in such a way that $\f{d}{d'}\to\la$, the symmetrization of the singular law  of  $M_1(\nu_n)+\cdots +M_{k_n}(\nu_n)$ (with $M_1(\nu_n),\ldots, M_{k_n}(\nu_n)$ independent and distributed according to $\mathbb{Q}_{d,d'}^{\nu_n}$) goes to $\nu_n^{\arc k_n}$. 

Moreover, if one fixes $d,d'$ and lets $n$ go to infinity, the distribution  $\ds\lf(\mathbb{Q}_{d,d'}^{\nu_n}\ri)^{*k_n}$ of $M_1(\nu_n)+\cdots+ M_{k_n}(\nu_n)$ converges weakly to a distribution   $\mathbb{P}_{d,d'}^{\mu}$ on the set of \dtpd matrices, whose Fourier transform is given by the following formula: for any  \dtpd matrix $A$ \begin{equation}\label{trFdeLd}
\E_{\mathbb{P}_{d,d'}^{\mu}}\lf(\exp\lf( i\Re \lf(\Tr A^*M\ri)\ri)\ri)=\exp\lf(\E\lf(d\tii\psi_\mu\lf(\Re\lf(<u,Av>\ri)\ri)\ri)\ri)
\end{equation}
where 
$\psi_\mu$ is the {\it L\'evy exponent} of $\mu$, i.e. the unique continuous function $f$ on $\R$ \st $f(0)=0$ and the Fourier transform of $\mu$ is $\exp\circ f$,  
$<.,.>$ is the canonical hermitian product of $\C^{d}$, 
and $u=(u_1,\ldots,u_d)$, $v=(v_1,\ldots,v_{d'})$ are independent random vectors, uniformly distributed on the unit sphere of respectively $\C^d$, $\C^{d'}$.
The proof of this weak convergence,  analogous to the one of \teo 3.1 of \cite{fbg.ID},  uses the polar decomposition of $d \tii d'$ matrices and the bi-unitarily invariance of the distributions $\mathbb{Q}_{d,d'}^{\nu_n}$. Note that for all $\mu,\nu$, $\mathbb{P}_{d,d'}^\mu*\mathbb{P}_{d,d'}^\nu=\mathbb{P}_{d,d'}^{\mu *\nu}$, and that when $\mu=N(0,1)$, $\mathbb{P}_{d,d'}^\mu$ is the distribution of a matrix $[M_{i,j}]$ with $\lf(\Re M_{i,j},\Im M_{i,j}\ri)_{\substack{1\leq i\leq d\\ 1\leq j\leq d'
}}$  i.i.d. random variables $N(0,\ff{2d'})$-distributed.

So the convergence of the symmetrization of the singular law of a \pmu-random matrix is the expression of the commutativity of the following diagram: $$\begin{array}{ccc}M_1(\nu_n)+\cdots+ M_{k_n}(\nu_n)&\;\ds \stackrel{n\to \infty}{-\!\!\!-\!\!\!-\!\!\!\longrightarrow}\; &\mathbb{P}_{d,d'}^{\mu}\\
\arrowvert&&\arrowvert\\
\substack{{\textrm{$d,d'$ go to $\infty$}}\\ d /d'\simeq \la}&&\substack{{\textrm{$d,d'$ go to $\infty$}}\\ d/d'\simeq \la}\\
\downarrow &&\downarrow\\
\substack{{\textrm{symmetrized}}\\ {\textrm{singular law:}}\\ \nu_n^{\arc k_n}}&\;\ds \stackrel{n\to \infty}{-\!\!\! -\!\!\!-\!\!\!\longrightarrow}\; &\substack{{\textrm{symmetrized}}\\ {\textrm{singular law:}}\\ \La_\la(\mu)}
\end{array}$$

To prove this result, we need a preliminary result about cumulants of $\arc$-infinitely divisible laws with compactly supported L\'evy measure. First, note that by lemma \ref{29.04.04.1},  such laws are compactly supported. Recall that free cumulants with ratio $ \la$ have been defined in the beginning of section \ref{6.2.05.1} by (\ref{PlayItSam}). For $\nu$  \pro measure $\nu$ whose moments of all orders are defined, let us denote by $c^*_n(\nu)$ ($n\geq 1$) its classical cumulants. Recall that classical cumulants  linearize the classical convolution and  are linked to the moments by the formula \begin{equation}\label{22.9.06.1}\ds\forall k\geq 1, \, m_k(\nu)=\sum_{\pi\in \Part (k)}\prod_{V\in \pi} {c}_{|V|}^*(\nu).\end{equation}
\begin{Th}\label{perig2}
Let $\mu$ be a symmetric \cid with compactly supported L\'evy measure, and let, for each integer $n$, $\mu_n$ be a \pro measure \st $\mu_n^{* n}=\mu$.
Then for each $k\geq 1$, the sequence $(n\tii  m_{2k}(\mu_n))_n$ tends to the $2k$-th classical cumulant $c_{2k}^*(\mu)$ of $\mu$, which is equal to $c_{2k}(\La_\la(\mu))$. 
\end{Th}

\begin{pr} By (\ref{22.9.06.1}), for all $n$, one has $$n\tii m_{2k}(\mu_n)=n\!\!\sum_{\pi\in\Part (2k)}\prod_{V\in \pi}\underbrace{{c}_{|V|}^*(\mu_n)}_{n^{-1}{c}_{|V|}^*(\mu)}=\sum_{\pi\in\Part (2k)}n^{1-|\pi |}{c}_{\pi}^*(\mu)={c}_{2k}^*(\mu)+o(1).$$
Let us denote $\nu_n=\mu_n^{\sst\boxplus_\la\ds n}$. By the line above, for all $k$, $m_{2k}(\mu_n)=O(n^{-1})$, so, by an easy induction on $k$ based on equation (\ref{PlayItSam}), one gets $c_{2k}(\mu_n)=O(n^{-1})$. Hence $c_{2k}(\nu_n)=O(1)$, and $m_{2k}(\nu_n)=O(1)$.  Moreover, by the equivalence (\ref{11.12.03.1}), the sequence $(\nu_n)$ converges weakly to $\La_\la(\mu)$. So the moments of $\nu_n$ tend to the moments of $\La_\la(\mu)$ (cf \cite{billingsley}). But thanks to (\ref{PlayItSam}),   \begin{eqnarray*}\ds n\tii m_{2k}(\mu_n)&=&n\sum_{\pi\in\NC' (2k)}\la^{e(\pi)}\prod_{V\in \pi}\underbrace{c_{|V|}(\mu_n)}_{n^{-1}c_{|V|}(\nu_n)}\\ & =&\sum_{\pi\in\NCp (2k)}\la^{e(\pi)}n^{1-|\pi |}\prod_{V\in \pi} {c}_{|V|}(\nu_n).\end{eqnarray*}  It has already been proved just above that the left hand term of the previous equation  tends to $c_{2k}^*(\mu)$, whereas the right hand term tends to $$\ds\sum_{\pi\in\NCp (2k)}\la^{e(\pi)}\delta_1^{|\pi |}\prod_{V\in \pi} {c}_{|V|}(\La_\la(\mu))={c}_{2k}(\La_\la(\mu)).$$It allows us to conclude.
\end{pr}

We will first prove the result when $\mu$ has a compactly supported L\'evy measure. We will work with a sequence $(d'_d)_d$ \st $1\leq d\leq d'_d$, and $d/d'_d$ tends to $\la\in (0,1]$ (even though the proof can be adapted to the case $\la=0$, we assume that $\la>0$ in order to simplify). To simplify notations, $d'$ will stand for $d'_d$.
\begin{propo}
Let $\mu$ be a symmetric \cid with compactly supported L\'evy measure. Then for all integers $k$,
\\
\\
(a) $\ds\lim_{d\to\infty}\E_{\pmum}\lf(m_k(\tilde{\mu}_{|M|})\ri)=m_k\lf(\La_\la(\mu)\ri).$
\\
\\
(b) The variance, under $\mathbb{P}_{d,d'}^\mu$, of $m_k(\tilde{\mu}_{|M|})$ tends to zero as $d$ goes to infinity. 
\end{propo}

\begin{pr} For an integer $n$, $[n]$ will denote $\{1,\ldots, n\}$, and $\NC(n)$ will denote the set of noncrosing partitions of $[n]$. Recall that $\NCp(n)$ denotes the set of noncrosing partitions of $[n]$ in which all blocks have even cardinality.

(a) First, for every complex $d\tii d'$ matrix $M$, for all  integer $k$, $m_k\lf(\tilde{\mu}_{|M|}\ri)$ is null if $k$ is odd and is equal to $\tr \lf(MM^*\ri)^{\f{k}{2}}$ ($\tr$ denotes normalized trace) if $k$ is even.  $\La_\la(\mu)$ is symmetric, so it suffices to show that for all   $k\in\n^*$,  $$\ds\lim_{d\to\infty}\E_{\pmum}\lf(\tr \lf(MM^*\ri)^k\ri)=m_{2k}\lf(\La_\la(\mu)\ri) .$$ 
Let, for $n\in\n^*$, $\mu_n$ be the \pro measure \st $\mu_n^{*n}=\mu$. Consider, for $d\geq 1$ and $n\geq 1$, $\lf(M_{d,n}^{(i)}\ri)_{1\leq i\leq n}$ i.i.d. random matrices with distribution \rdmunp By definition, for every $d\geq 1$, the sum of the $M_{d,n}^{(i)}$'s ($i=1\ldots n$) converges in distribution to \pmu when  $n$ goes to $\infty$.
\\
We know, by \teo \ref{perig2} that for all $k\in\n^*$, the sequence $(n\tii m_k(\mu_n))_n$ is bounded, and so (see \cite{billingsley})  for all $k,d\in\n^*$, \begin{equation}\label{perig42}
\E_{\pmum}\lf(m_{2k}\lf(\tilde{\mu}_{|M|}\ri)\ri)=\lim_{n\to\infty}\E\lf(\tr\lf(\lf(\sum_{i=1}^nM_{d,n}^{(i)}\ri)\lf(\sum_{i=1}^nM_{d,n}^{(i)*}\ri)\ri)^k\ri) .
\end{equation}
 {\bf Let us fix } $\mathbf{k\in \n^\ast}$.
\\
We are going to use the formula (\ref{perig42}).
\\
Let, for $d,n\geq 1$, \begin{equation*}
b_{d,n}=\E\lf(\tr\lf(\lf(\sum_{i=1}^nM_{d,n}^{(i)}\ri)\lf(\sum_{i=1}^nM_{d,n}^{(i)*}\ri)\ri)^k\ri).\end{equation*}
From now on, we do not write anymore the index $d$ in  $M_{d,n}^{(i)}$. We denote, for $l,n$ non-negative integers, by $A_n^l$  the number of one-to-one maps from $[l]$ to $[n]$, i.e. $n(n-1)\cdots(n-l+1)$. For a partition $\pi$  of $[2k]$, 
for $1\leq l\leq 2k$, we denote by $\pi(l)$ the  index of the class
of $l$,  after having ordered the classes according to the order of their first element (for example, $\pi(1)=1$;  $\pi(2)=1$ if $1\stackrel{\pi}{\sim}2$ and  $\pi(2)=2$ if $1\stackrel{\pi}{\nsim}2$). Then we have 
\begin{eqnarray*}
b_{d,n}&=& \tr \lf(\E\lf(\sum_{\scriptstyle{f\in\{1,\ldots,n\}^{2k}}}M_{n}^{(f(1))}M_{n}^{(f(2))*}\cdots M_{n}^{(f(2k))*}\ri) \ri)\\
&=& \tr \lf( \E\lf(\sum_{\scriptstyle\pi\in \Part(2k)}\ds A_{n}^{|\pi|} M_{n}^{(\pi(1))}M_{n}^{(\pi(2))*}M_{n}^{(\pi(3))}\cdots M_{n}^{(\pi(2k))*}\ri)\ri).\end{eqnarray*} 
\[\textrm{But }\quad\quad\E\lf(\underbrace{M_n^{(1)*}M_n^{(1)}\cdots M_n^{(1)*}}_{\textrm{$2l+1$ alterned factors}}\ri)=\E\lf(\underbrace{M_n^{(1)}M_n^{(1)*}\cdots M_n^{(1)}}_{\textrm{$2l+1$ alterned factors}}\ri)=0,\]
\[\E\lf(\underbrace{M_n^{(1)*}M_n^{(1)}\cdots M_n^{(1)}}_{\textrm{$2l$ alterned factors}}\ri)=\f{d}{d'}m_{2l}(\mu_n)I_{d'},\]
\[\E\lf(\underbrace{M_n^{(1)}M_n^{(1)*}\cdots M_n^{(1)*}}_{\textrm{$2l$ alterned factors}}\ri)=m_{2l}(\mu_n)I_{d}.\]
So, using many times the fact that a partition $\pi$ of a finite totally ordered set is noncrossing if and only if one of its class $V$ is an interval and $\pi\backslash \{V\}$ is noncrossing (page 3 of \cite{spei98}) and integrating  successively with respect to the different independent random matrices, one has  
$$\ds\pi\in \NCp(2k)\Rightarrow \tr\E \lf(M_{n}^{(\pi(1))}M_{n}^{(\pi(2))*}\cdots M_{n}^{(\pi(2k))*}\ri)=\ds\lf(\f{d}{d'}\ri)^{e(\pi)}\prod_{B\in \pi} m_{|B|}(\mu_n),$$ $$
\pi\in \NC(2k)\backslash \NCp(2k)\Rightarrow \tr\E \lf(M_{n}^{(\pi(1))}M_{n}^{(\pi(2))*}\cdots M_{n}^{(\pi(2k))*}\ri)=0.$$

But $A_n^{|\pi|}\sim n^{|\pi|}$, so, by the preceding \teov the limit, when $n$ goes to infinity, of $$\tr \lf( \E\lf(\sum_{\scriptstyle\pi\in \Part(2k)}\ds A_{n}^{|\pi|} M_{n}^{(\pi(1))}M_{n}^{(\pi(2))*}\cdots M_{n}^{(\pi(2k))*}\ri)\ri)$$ is $$\ds\sum_{\pi\in \NCp(2k)}\lf(\f{d}{d'}\ri)^{e(\pi)}\prod_{V\in \pi}c_{|V|}(\La_\la(\mu)),$$ which tends, when $d$ goes to infinity,  to $$\ds\sum_{\pi\in \NCp(2k)}\la^{e(\pi)}\prod_{V\in \pi}c_{|V|}(\La_\la(\mu)),$$which is equal to $m_{2k}\lf(\La_\la(\mu)\ri)$ by (\ref{PlayItSam}).

So it suffices to prove that $$\ds b_{d,n}':= \sum_{\substack{\pi\in \Part (2k)\\ \pi\notin \NC(2k)}} A_{n}^{|\pi|}\tr \E\lf\{  M_{n}^{(\pi(1))}M_{n}^{(\pi(2))*}M_{n}^{(\pi(3))}\cdots M_{n}^{(\pi(2k))*}\ri\}$$ vanishes when $n$, and then $d$, go to infinity. Let us expand the trace: $b'_{n,d}$ is equal to$$\ds\frac{1}{d}
\sum_{
\substack{
\pi\in\Part (2k)\\ \pi\notin\NC (2k)
}}
\sum_{
\sst\substack{j\in [d']^{2k}\\ \forall r \textrm{ odd, }j_r\leq d}\ds
}A_n^{|\pi|}\E\lf\{
\lf[M_{n}^{(\pi(1))}\ri]_{j_1,j_2}\lf[M_{n}^{(\pi(2))*}\ri]_{j_2,j_3}\cdots \lf[M_{n}^{(\pi(2k))*}\ri]_{j_{2k},j_1}
\ri\}.
$$
Using the fact that     $\lf(M_{d,n}^{(i)}\ri)_{1\leq i\leq n}$ are independent copies of a matrix with distribution $\mathbb{Q}_{d,d'}^{\mu_n}$, we deduce  (with the notation $j_{2k+1}=j_1$)\begin{eqnarray*}
b'_{d,n}&=&\frac{1}{d}
\sum_{\substack{\pi\in\Part (2k)\\ \pi\notin\NC (2k)}}
A_{n}^{|\pi|}\sum_{
\sst\substack{
j\in [d']^{2k}\\ \forall r \textrm{ odd, }j_r\leq d
}}
\prod_{B\in \pi}
\E_{\mathbb{Q}_{d,d'}^{\mu_n}}\lf\{
\prod_{\substack{r\in B\\ \textrm{$r$ odd}}} M_{j_r,j_{r+1}}\prod_{\substack{r\in B\\ \textrm{$r$ even}}} M_{j_r,j_{r+1}}^* \ri\}
\\
&=&\frac{1}{d}
\sum_{
\substack{
\pi\in\Part (2k)\\ \pi\notin\NC (2k)}}\f{A_{n}^{|\pi|}}{n^{|\pi|}}
\sum_{
\scriptstyle\substack{j\in [d']^{2k}\\ \forall r \textrm{ odd, }j_r\leq d }\ds
}\\ &&\prod_{
\scriptstyle B\in \pi\ds}
n\E\lf\{\sum_{\scriptstyle l\in[d]^B\ds}
\prod_{\scriptstyle \substack{r\in B\\ \textrm{$r$ odd}}\ds} 
u_{j_r,l_r} v_{l_r,j_{r+1}}X_{n,l_r} \prod_{\substack{r\in B\\ \textrm{$r$ even}}} \bar{v}_{l_r,j_r}X_{n,l_r}  \bar{u}_{j_{r+1},l_r}\ri\},
\end{eqnarray*}
where $U,V ,X_{n,1},\ldots,X_{n,d}$ are independent, with respective distribution  the Haar measure on the group of $d\tii d$ unitary matrices, the Haar measure on the group of $d'\tii d'$ unitary matrices, and $\mu_n$.
 
For all $B\subset [2k]$, for all $j\in [d']^{2k}$ \st for all $r$ odd, $j_r\leq d$, summing over the partition generated by $l$, one has \begin{eqnarray*}&\ds n\E\lf\{\sum_{\scriptstyle l\in[d]^B\ds}
\prod_{\scriptstyle \substack{r\in B\\ \textrm{$r$ odd}}\ds} 
u_{j_r,l_r} v_{l_r,j_{r+1}}X_{n,l_r} \prod_{\substack{r\in B\\ \textrm{$r$ even}}} \bar{v}_{l_r,j_r}X_{n,l_r}  \bar{u}_{j_{r+1},l_r}\ri\}&\\
=&\ds\sum_{\sigma\in\Part(B)}\sum_{\scriptstyle \substack{l\in[d]^\sigma}\ds}n^{1-|\sigma|}\lf[\prod_{W\in \sigma} n m_{|W|}(\mu_n)\ri]\E\lf\{
\prod_{\scriptstyle \substack{r\in B\\ \textrm{$r$ odd}}\ds} 
u_{j_r,l_r} v_{l_r,j_{r+1}} \prod_{\substack{r\in B\\ \textrm{$r$ even}}} \bar{v}_{l_r,j_r}  \bar{u}_{j_{r+1},l_r}\ri\}. &\end{eqnarray*}
The measures $\mu_n$ are symmetric, so all their moments of odd order are null. Hence, according to theorem \ref{perig2}, the quantity of the previous equation tends, as $n$ goes to infinity, to $$\begin{cases}\ds\sum_{l=1}^dc_{|B|}^*(\mu) \E\lf\{
\prod_{\scriptstyle \substack{r\in B\\ \textrm{$r$ odd}}\ds} 
u_{j_r,l} v_{l,j_{r+1}} \prod_{\substack{r\in B\\ \textrm{$r$ even}}} \bar{v}_{l,j_r}  \bar{u}_{j_{r+1},l}\ri\}&\textrm{if $|B|$ is even,}\\  0&\textrm{otherwise.}\end{cases}
$$ Note that equation (\ref{22.9.06.1}) implies, by an easy induction,  that the classical cumulants with of even order of a symmetric \pro measure are null. It allows us to avoid considering two cases on the parity of $|B|$ in the previous limit, and to claim that the first formula is valid whenever $|B|$ is odd.  
This limit  is equal, by invariance of uniform distributions on unitary groups by permutation of rows and columns, to $$\ds d \tii c_{|B|}^*(\mu) \E\lf\{
\prod_{\scriptstyle \substack{r\in B\\ \textrm{$r$ odd}}\ds} 
u_{j_r} v_{j_{r+1}} \prod_{\substack{r\in B\\ \textrm{$r$ even}}} \bar{v}_{j_r}  \bar{u}_{j_{r+1}}\ri\},$$ where $u,v$ are independent uniform random vectors on the unit spheres of $\C^d,\C^{d'}$.

So the limit, when $n$ goes to infinity, of $b'_{d,n}$ is   $$\ds\frac{1}{d}
\sum_{
\substack{
\pi\in\Part (2k)\\ \pi\notin\NC (2k)}}
\sum_{
\scriptstyle\substack{j\in [d']^{2k}\\ \forall r \textrm{ odd, }j_r\leq d }\ds
}\prod_{
\scriptstyle B\in \pi\ds}d\tii c_{|B|}^*(\mu) \E\lf\{
\prod_{\scriptstyle \substack{r\in B\\ \textrm{$r$ odd}}\ds} 
u_{j_r} v_{j_{r+1}} \prod_{\substack{r\in B\\ \textrm{$r$ even}}} \bar{v}_{j_r}  \bar{u}_{j_{r+1}}\ri\},$$ the absolute value of which is less or equal, by invariance of the distributions of $u$ and $v$ under permutation of coordinates, to   $$\ds\frac{1}{d}
\sum_{
\substack{
\pi,\tau\in\Part (2k)\\ \pi\notin\NC (2k)}}
{d'}^{|\pi|+|\tau|}\prod_{
\scriptstyle B\in \pi\ds}\lf|c_{|B|}^*(\mu) \E\lf(\prod_{\underset{\scriptstyle{r\:odd}}{\scriptstyle{r\in B}}}u_{\tau(r)}v_{\tau(r+1)}\prod_{\underset{\scriptstyle{r\:even}}{\scriptstyle{r\in B}}}\overline{u}_{\tau(r+1)}\overline{v}_{\tau(r)}\ri)\ri|.$$
Moreover, by invariance of the distribution of $u$ under the action of unitary diagonal matrices, for every pair $(\pi,\tau)$ of partitions of  $[2k]$, if 
$$\prod_{B\in \pi}\E\lf(\prod_{\underset{\scriptstyle{r\:odd}}{\scriptstyle{r\in B}}}u_{\tau(r)}v_{\tau(r+1)}\prod_{\underset{\scriptstyle{r\:even}}{\scriptstyle{r\in B}}}\overline{u}_{\tau(r+1)}\overline{v}_{\tau(r)}\ri)$$ 
is non zero, then for every class $B$ of $\pi$, there exists $\phi$, permutation of $B$, which maps odd numbers to even numbers and vice versa, \st for all  $r\in B$, $\tau(r)=\tau(\phi(r)+1)$. It  implies, by lemma 4.4 of \cite{fbg.ID}, that  $|\tau|+|\pi|\leq 2k$. So one has, using the H\"older inequality and equation (4.2.11) of \cite{hiai},  $$\ds\lim_{n\to\infty}b'_{n,d}=O(d^{-1}),$$
which closes the proof of (a).

One notes that the proof of $(a)$ is a very closed adaptation of the proof of Proposition 4.1 of \cite{fbg.ID}, by adaptation of the arguments to the context of  non hermitian and non square matrices. Using again the same adaptation, the proof of $(b)$ is along the same lines as the proof of Proposition 5.1 of \cite{fbg.ID}. 
\end{pr}

To conclude this section, we have to state its main theorem. Recall that the {\it convergence in probability} of a sequence $X_n$ of random variables in a metrizable topological space $ \mc{X}$ to a constant $l\in \mc{X}$ is the convergence of the probability of the event $\{d(X_n,l)<\eps\}$ to $1$ for all positive $\eps$, where $d$ is any distance which defines the topology (it does not depend on the choice of such a distance). In the following theorem, we shall refer to the convergence in \pro in the set of \pro measures on the real line, endowed with the metrizable topology of weak convergence.  The proof of the \teo  is similar to the one of Theorem 7.6 of \cite{fbg.ID}, based on the previous proposition and on  an approximation by compound Poisson laws. The only modification is to work with products of the type $MM^*$ rather than $M^*M$.  Recall that $d'$ is in fact $d'_d$ (i.e. $d'$ depends on $d$) and that $d/d'$ tends to $\la$ as $d$ tends to infinity.

\begin{Th}\label{CVPROBAberco2}
Let $\mu$ be a symmetric \cidp Let, for $d\geq 1$, $M_d$ be a random matrix with distribution  $\mathbb{P}_{d,d'}^{\mu}$.
\\
Then as $d$ goes to infinity,  the symmetrization $\tilde{\mu}_{|M_d|}$ of the singular law of $M_d$  converges in probability to $\La_\la(\mu)$. 
\end{Th}

\begin{rmq}In the case where $\mu$ is a normal law, one recovers the well known result about asymptotic spectral distribution of Wishart random matrices.\end{rmq}

\section{Rectangular symmetric Poisson distributions as limits of sums of rank-one matrices}\label{france.q=chinois.a.Belleville}
 The symmetric Poisson distribution with parameter $c>0$ has been introduced in section \ref{bon.vieux.dede.20.09.06}. The  free analogues of unsymmetric Poisson distributions are Marchenko-Pastur laws. But as we said it in section \ref{bon.vieux.dede.20.09.06}, unless $\la=0$, the computations for symmetric Poisson laws are harder than for the unsymmetric ones (even when $\la=1$, the densities have not been expressed). Nevertheless, we have the following  characterization of the rectangular analogues of symmetric Poisson distributions.

\begin{propo}Consider $\la\in (0,1]$, and $c>0$. Then the image, by the Bercovici-Pata bijection with ratio $\la$, of the symmetric Poisson distribution with parameter $c$ is the limit, for convergence in probability, of the symmetrization of the singular law of the random matrix $$M(d,d',d''):=\ds\sum_{k=1}^{d''}u_d(k)v_{d'}(k)^*$$ 
when \begin{equation}\label{6.2.05.2}d\to\infty,\quad \f{d}{d'}\to \la,\quad \f{d''}{d}\to c,\end{equation} 
where $u_d(k), v_{d'}(k)$ ($k\geq 1$) are independent uniform random vectors of the unit spheres of $\C^d,\C^{d'}$ (considered as column matrices). 
\end{propo}

\begin{rmq}Note that when $\la=1$, the image, by the Bercovici-Pata bijection with ratio $\la$, of the symmetric Poisson distribution with parameter $c$ is $\tau_+\bxp \tau_-$, where $\tau_+$ is the  Marchenko-Pastur distribution with parameter $c/\! 2$, and $\tau_-$ is the push-forward, by the function $t\to -t$, of $\tau_+$. \end{rmq}

\begin{pr}Let us denote by $\mu$ the  symmetric Poisson distribution with parameter $c$, and by $\sigma$ the push-forward by $t\to t^2$ of its  image by the Bercovici-Pata bijection with ratio $\la$. As explained in the beginning of the proof of the last theorem of the section called "the rectangular  $R$-transform" of \cite{fbg.AOP.rect}, it suffices to prove that  for each $\eps>0$, the \pro of the event $$\lf\{\underset{\Im z\geq 1}{\sup}\lf|\ff{d}\Tr\mf{R}_z\lf(A(d,d',d'') \ri)-G_{\sigma}(z)\ri|>\eps\ri\}$$ 
tends to zero as $d, d', d''$ tend to infinity as in (\ref{6.2.05.2}), where $$A(d,d',d'')=M(d,d',d'')M(d,d',d'')^*, $$ and for $M$ hermitian matrix and $z$ complex non real number, $\mf{R}_z\lf(M\ri)=(z-M)^{-1}$. 

Fix $\eps >0$.
It can easily be seen that $\mathbb{P}_{d,d'}^{\mu}$ is the distribution of  $$N(d,d'):=\ds\sum_{k=1}^{X(cd)}u_d(k)v_{d'}(k)^*,$$where $X(cd)$ is a random variable distributed according to an (unsymmetric) Poisson distribution with parameter $cd$ and $X(cd)$ is independent of the $u_d(k)$'s and of the $v_{d'}(k)$'s. Thus by the previous theorem, with the notation $ B(d,d')=N(d,d')N(d,d')^*,$
 the probability of the event  $$\lf\{\underset{\Im z\geq 1}{\sup}\lf|\ff{d}\Tr\mf{R}_z\lf(B(d,d') \ri)-G_{\sigma}(z)\ri|>\eps\ri\}$$ tends to zero. Thus it suffices to prove that the probability of the event  $$\lf\{\underset{\Im z\geq 1}{\sup}\lf|\ff{d}\Tr\lf(\mf{R}_z\lf(B(d,d') \ri)- \mf{R}_z\lf(A(d,d',d'')\ri)\ri)\ri|>\eps\ri\}$$ tends to zero. But for all hermitian $d\tii d$ matrices $A,B$, for all $z$ \st $\Im z\geq 1$,  we have $$\mf{R}_z(B)-\mf{R}_z(A)=-\mf{R}_z(B)(B-A)\mf{R}_z(A),$$ whose normalized trace is not more than its norm times its rank divided by $d$. Moreover, $||\mf{R}_z(B)-\mf{R}_z(A)||\leq 2$, and the rank is not more than the one of $B-A$. So it suffices to prove that $$\ff{d}\rg\lf(B(d,d')-A(d,d',d'')\ri)$$ converges in \pro to zero. $B(d,d')-A(d,d',d'')$ can be put in the form $$(\ldots)(N(d,d')-M(d,d',d''))^*+(N(d,d')-M(d,d',d''))(\ldots),$$ so $$\ff{d}\rg\lf(B(d,d')-A(d,d',d'')\ri)\leq \f{2}{d}\rg(N(d,d')-M(d,d',d''))\leq \f{2}{d}|X(cd)-d''|,$$ which converges in \pro to zero, by the weak law of great numbers.
\end{pr}

Florent Benaych-Georges\\ DMA, \'Ecole Normale Sup\'erieure,\\ 45 rue d'Ulm, 75230 Paris Cedex 05, France\\   e-mail: benaych@dma.ens.fr\\
  http://www.dma.ens.fr/$\sim$benaych

\begin{thebibliography}{99}
\bibitem[A61]{akhi} Akhiezer, N.I. \emph{The classical moment problem}, Moscou, 1961
\bibitem[B-NT02]{steen2} Barndorff-Nielsen, O.E., Thorbj\o rnsen, S.  \emph{Selfdecomposability and Levy processes in free probability}, Bernoulli 8(3) (2002), 323-366. 
\bibitem[B-G04]{fbg.ID} Benaych-Georges, F. \emph{Classical and free infinitely divisible distributions and random matrices} Annals of Probability, Vol. 33, No. 3, p. 1134-1170 (2005), available on the web page of the author
\bibitem[B-G1]{fbg.AOP.rect} Benaych-Georges, F. \emph{Rectangular random matrices. Related convolution} submitted, available on the web page of the author\bibitem[B-G2]{fbg.free.amalg} Benaych-Georges, F. \emph{Rectangular random matrices, related free entropy and free Fisher's information} submitted, available on the web page of the author
\bibitem[BPB99]{appenice} Bercovici, H., Pata, V., with an appendix by Biane, P. \emph{Stable laws and domains of attraction in free
probability theory} Annals of Mathematics, 149 (1999) 1023-1060
\bibitem[BV93]{defconv} Bercovici, H., Voiculescu, D. \emph{Free convolution of measures with unbounded supports} Indiana Univ.
Math. J. 42 (1993) 733-773
\bibitem[B68]{billingsley} Billingsley, P. \emph{Convergence of probability measures} Wiley, 1968
\bibitem[C-D04]{cabduv.ID} Cabanal-Duvillard T. \emph{A matrix representation of the Bercovici-Pata bijection} Electron. J. Probab.  10  (2005), no. 18, 632--66
\bibitem[D74]{donog} Donoghue, W. \emph{Monotone matrix functions and analytic continuation}, Springer, New-York, 1974
\bibitem[F66]{feller2} Feller, W. \emph{An introduction to probability theory and its applications}, volume II, second edition, New York London Sydney : J. Wiley, 1966
\bibitem[GK54]{gne} Gnedenko, V., Kolmogorov, A.N. \emph{Limit distributions for sums of independent random variables}
Adisson-Wesley Publ. Co., Cambridge, Mass., 1954
\bibitem[HP00]{hiai} Hiai, F., Petz, D. \emph{The semicircle law, free random variables, and entropy} Amer.
Math. Soc., Mathematical Surveys and Monographs Volume 77, 2000
 \bibitem[N74]{Nelson} Nelson, E.  \emph{Notes on non-commutative integration} J.
Functional Analysis 15 (1974), 103--116
\bibitem[OP97]{orpetz} Oravecz, F., Petz, D. \emph{On the eigenvalue distribution of some symmetric random matrices} Acta Sci. Math. (Szeged) 63 (1997), no. 3-4, 383--395.
\bibitem[P97]{petrov} Petrov, V.V. \emph{Limit theorems of probability theory} Oxford Studies in Probability, 4, 1995
\bibitem[S98]{spei98} Speicher, R. \emph{Combinatorial theory of the free product with amalgamation and operator-valued free probability theory}, Mem. Amer. Math. Soc. 132 (1998), no. 627
\bibitem[VDN91]{voic1} Voiculescu, D.V., Dykema, K., Nica, A. \emph{Free random variables} CRM Monograghs Series No.1, Amer.
Math. Soc., Providence, RI, 1992 
\end{thebibliography}
\end{document}